\documentclass[amsthm,secthm,11pt,preprint]{elsart}
\usepackage[english]{babel}
\usepackage[latin1]{inputenc}
\usepackage{float} 
\usepackage[english, algosection, algoruled, noline]{algorithm2e}
\usepackage{latexsym}
\usepackage{amsmath}
\usepackage{amssymb}
\usepackage{graphicx}
\usepackage{amsfonts}
\usepackage{color}
\usepackage{yjsco}
\usepackage{natbib}
\usepackage[T1]{fontenc}

\usepackage{pdflscape}

\newtheorem{example}[thm]{Example}

\newcommand{\R}{\mathbb{R}}

\newcommand{\C}{\mathbb{C}}
\newcommand{\N}{\mathbb{N}}

\newcommand{\Ac}{\mathcal{A}} 
\newcommand{\Lc}{\mathcal{L}} 
\newcommand{\Qc}{\mathcal{Q}} 
\newcommand{\Sc}{\mathcal{S}}
\newcommand{\Pc}{\mathcal{P}} 
\newcommand{\Ic}{\mathcal{I}}

\newcommand{\Vc}{\mathcal{V}} 
 
\newcommand{\BB}{\mathcal{B}}

\newcommand{\Mon}{\mathcal{M}}
\newcommand{\Span}[1]{\<{#1}\>}

\newcommand{\gb}{\mathbf g}
\newcommand{\xx}{\mathbf x}
\newcommand{\uu}{\mathbf u}
\newcommand{\Fc}{\mathcal{F}}

\newcommand{\SpanD}[2]{\<#1\,|\,{ #2}\>}

\newcommand{\piFB}{\pi_{F,\BB}}
\newcommand{\rank}{\mathrm{rank}}
\newcommand{\ev}{\ensuremath{\boldsymbol{1}}}

\def\<{\langle}
\def\>{\rangle}

\newcommand{\nobracket}{}

\newcommand{\plusassign}{+\!\!=}

\newcommand{\tmstrong}[1]{\textbf{#1}}

\newcommand{\udots}{{\mathinner{\mskip1mu\raise1pt\vbox{\kern7pt\hbox{.}}\mskip2mu\raise4pt\hbox{.}\mskip2mu\raise7pt\hbox{.}\mskip1mu}}}

\newenvironment{Proof}{\noindent\textbf{Proof.\ }}{\hspace*{\fill}$\Box$\medskip}

\newenvironment{changemargin}[2]{%
\begin{list}{}{%
\setlength{\topsep}{0pt}%
\setlength{\leftmargin}{#1}%
\setlength{\rightmargin}{#2}%
\setlength{\listparindent}{\parindent}%
\setlength{\itemindent}{\parindent}%
\setlength{\parsep}{\parskip}%
}%
\item[]}{\end{list}} 

\begin{document}

\begin{frontmatter}
%\Hacemos una portada
\title{Border basis relaxation for polynomial optimization}

\author[inria]{Marta Abril Bucero}
\ead{Marta.Abril\_Bucero@inria.fr},
\author[inria]{Bernard Mourrain} 
\ead{Bernard.Mourrain@inria.fr}
\address[inria]{EPI GALAAD, INRIA M\'editerran\'ee, 2004 Route des Lucioles,
         BP 93, 06902 Valbonne, France}

%\thispagestyle{empty}
%\newpage

%\Comenzamos con un resumen
\begin{abstract}
  A relaxation method based on border basis reduction which improves
  the efficiency of Lasserre's approach is proposed to compute the
  infimum of a polynomial function on a basic closed semi-algebraic
  set. A new stopping criterion is given to detect when the relaxation
  sequence reaches the infimum, using a sparse flat extension
  criterion. We also provide a new algorithm to reconstruct a finite
  sum of weighted Dirac measures from a truncated sequence of moments,
  which can be applied to other sparse reconstruction problems. As an
  application, we obtain a new algorithm to compute zero-dimensional
  minimizer ideals and the minimizer points or
  zero-dimensional G-radical ideals. Experiments show the impact of
  this new method on significant benchmarks.
\end{abstract}

\end{frontmatter}

%\newpage
%\Introducimos un \'indice y escribimos los contenidos de cada secci\'on usando la orden \section
%\tableofcontents
%\thispagestyle{empty}
%\newpage
\setcounter{page}{1}

\section{Introduction}

Computing the global infimum of a polynomial function $f$ on a
semi-algebraic set is a difficult but important problem, with many
applications. A relaxation approach was proposed in
\citep{Las01} (see also \citep{Par03}, \citep{Shor87}) which approximates
this problem by a sequence of finite dimensional convex optimization problems.
These optimization problems can be formulated in terms of
linear matrix inequalities on moment matrices associated to the set of
monomials of degree $\leq t\in \N$ for increasing values of $t$. They
can be solved by Semi-Definite Programming (SDP) techniques.
The sequence of minima converges to the actual infimum $f^{*}$ of the
function under some hypotheses \citep{Las01}. In some cases, the
sequence even reaches the infimum in a finite number of steps
\citep{Lau07, NDS, Marshall, DNP, Ha-Pham:10, Nie11}.
This approach has proved to be particularly fruitful in many problems \citep{Lasserre:book}.
In contrast with numerical methods such as gradient descent methods,
which converge to a local extremum but with no guaranty for the global
solution, this relaxation approach can provide certificates for the infimum
value $f^{*}$ in terms of sums of squares representations. 

From an algorithmic and computational perspective, however some issues need
to be considered. 

The size of the SDP problems to be solved is a bottleneck of the
method. This size is related to the number of monomials of degree $\leq t$
and increases exponentially with the number of variables and the
degree $t$. 
Many SDP solvers are based on interior point methods which provide an
approximation of the optimal moment sequence within a given precision
in a polynomial time: namely
$\mathcal{O} ((p\, s^{3.5}+ c\,p^{2} s^{2.5}+ c\, p^{3}s^{0.5})\log(\epsilon^{-1}))$ arithmetic operations
where 
$\epsilon>0$ is the precision of the approximation,
$s$ is the size of the moment matrices, 
$p$ is the number of parameters (usually of the order $s^{2}$)
and $c$ is the number of constraints
\citep{NesterovNemirovski94}.
Thus reducing the size $s$ or the number of parameters $p$ can 
significantly improve the performance of these relaxation methods. Some recent
works address this issue, using symmetries (see
e.g. \citep{Riener:2013:ESS:2448396.2448402}) or 
polynomial reduction (see e.g. \citep{lasserre:hal-00651759}).
In this paper, we extend this latter approach.

While determining the infimum value of a polynomial function on a semi-algebraic
set is important, computing the minimizer points, is also critical in many applications.
Determining when and how these minimizer points can be computed from the relaxation
sequence is a problem that has been addressed, for instance in
\citep{HeLa05, Nie2012} using full moment matrices. 
This approach has been used for solving polynomial equations
\citep{Lau07, LLR07, LLR08b, Lasserre:book}.

The optimization 
problem can be reformulated as solving polynomial equations related to the (minimal) critical
value of the polynomial $f$ on a semi-algebraic set. 
Polynomial solvers based, for instance, on Gr\"obner basis
or border basis computation can then be used to recover the real critical points from
the complex solutions of (zero-dimensional) polynomial systems 
(see e.g. \citep{Parrilo03minimizingpolynomial, ElDin:2008:CGO:1390768.1390781,Greuet:2011:DRI:1993886.1993910}).
This type of methods relies entirely on polynomial algebra and
univariate root finding. So far, there is no clear comparison
of these elimination methods and the relaxation approaches.

\noindent{}\textbf{Contributions.}
We propose a new method which combines Lasserre's SDP relaxation approach
with polynomial algebra, in order to increase the
efficiency of the optimization algorithm. Border basis computations
are considered for their numerical stability
\citep{Mourrain2005,BMPhT08}. In principle, any graded normal form
technique could be used here.

A new stopping criterion is given to detect when the relaxation
sequence reaches the infimum, using a flat extension criterion from \citep{MoLa2008}. 
We also provide a new algorithm to reconstruct a finite sum of weighted Dirac
measures from a truncated sequence of moments. 
This reconstruction method can be used in other problems such
as tensor decomposition \citep{BracCMT09:laa} and multivariate sparse
interpolation \citep{Giesbrecht:2009:SSI:1539056.1539741}.

As shown in \citep{Cert,NDS,DNP,Marshall,Nie11,Ha-Pham:10},
an exact SDP relaxation can be constructed for ``well-posed'' optimization problems.
As an application, we obtain a new algorithm to compute zero-dimensional minimizer
ideals and the minimizer points, or zero-dimensional G-radicals.
Experiments show the impact of this new method compared to the previous
relaxation constructions.

\noindent{}\textbf{Content.}
The paper is organized as follows.  Section 2 describes the
minimization problem and includes a running example to explain
the different steps of our method. In Section 3-5, we describe the
ingredients of the main algorithm, which is described in Section 7.
In section 3, we describe the SDP relaxation hierarchies (full moment
matrices and border basis).
In Section 4, we tackle the sub-problem of how to compute the optimal
linear form through the solution of a SDP problem.  In Section 5, we
tackle the sub-problem of how to verify that we have found the infimum,
checking the flat extension property using orthogonal polynomials.  In
Section 6, we tackle the sub-problem of how to compute the minimizer
points using multiplication matrices.  Section 7 gives  a description
of the complete minimization algorithm. Section 8 analyses cases for
which an exact relaxation can be constructed.  Section 9
concludes experimentation.

\section{Minimization problem}
Let $f \in \R[\xx]$ be a polynomial function with real coefficients and let $\gb=\{g_1^{0},
\ldots, g_{n_{1}}^{0}$; $g_1^{+}, \ldots$, $g_{n_{2}}^{+}\}\in \R[\xx]$ be a set of constraints which is the union of a finite subset $\gb^{0}=\{g_1^{0},
\ldots, g_{n_{1}}^{0}\}$ of polynomials corresponding to the equality constraints and 
a finite subset $\gb^{+}=\{g_1^{+}, \ldots, g_{n_{2}}^{+}\}$ corresponding to
the non-negativity constraints.
The basic semi-algebraic set defined by the constraints $\gb$ will be denoted
$S:=\Sc (\gb) =\{ \xx \in \R^n \mid
g_1^{0}(\xx)=\cdots=g_{n_1}^{0}(\xx)=0, g_1^{+}(\xx)\ge  0,..., g_{n_2}^{+}(\xx)\ge 0
\}$.
We assume that $S\neq \emptyset$ and that $f$ is bounded by below on
$S$ (i.e. $\inf_{\xx \in S} f (\xx)>-\infty$).
%Let $f, g_1^{0}, \ldots, g_{n_{1}}^{0}$, $g_1^{+}, \ldots, g_{n_{2}}^{+}
%\in \R[\xx]$ be polynomial functions. 
The minimization problem that
we consider throughout the paper is the following: compute 
\begin{eqnarray}\label{problem1}
 \inf_{\xx \in \R^n} f(\xx) \\ 
 s.t. \ g_1^{0}(\xx)=\cdots=g_{n_1}^{0}(\xx)=0 \nonumber \\
     g_1^{+}(\xx)\ge 0,...,g_{n_2}^{+}(\xx)\ge 0  \nonumber
\end{eqnarray}
% More precisely, the objectives of the method we are going to describe
% are to compute the infimum value when
% $f$ is bounded by below and the minimizers points.

% A set of constraints $\gb=\{g_1^{0},
% \ldots, g_{n_{1}}^{0}$; $g_1^{+}, \ldots$, $g_{n_{2}}^{+}\}\in \R[\xx]$ is 
% the union of a finite subset $\gb^{0}=\{g_1^{0},
% \ldots, g_{n_{1}}^{0}\}$ of polynomials corresponding to the equality constraints and 
% a finite subset $\gb^{+}=\{g_1^{+}, \ldots, g_{n_{2}}^{+}\}$ corresponding to
% the non-negativity constraints.
% We say that $\gb \subset \gb'$ if $\gb^{0}\subset \gb'^{0}$ and $\gb^{+}\subset \gb'^{+}$.

When $n_1=n_2=0$, there is no constraint and $S=\R^{n}$. In this
case, we are considering a global unconstrained minimization problem.

The points $\xx^{*}\in \R^{n}$ which satisfy $f 
(\xx^{*}) = \inf_{\xx \in S} f (\xx)$ are called the
{\em minimizer points} of $f$ on $S$. 
The set of minimizer points is denoted
 $V_{min} = \{\xx^*\in S\ s.t \ f(\xx^*)= \inf_{\xx \in S} f (\xx) \}$.
The ideal of $\R[\xx]$ defining the set $V_{min}$ is denoted $I_{min}$.
The value $f^*=\inf_{\xx \in S} f (\xx)$ is called the {\em minimum}
of $f$ on $S$, when the set of minimizers is not empty.

If the set of minimizer points is not empty, we say that the {\em  minimization problem is feasible}.
The minimization problem is not feasible means that $V_{min}=\emptyset$ and $I_{min}=\R[\xx]$.

We present a running problem to explain the different steps of our method to solve Problem \eqref{problem1}.
\begin{example}
\begin{eqnarray*}
 \inf_{\xx \in \R^2} \ f(x,y)=(x-1)^2(x-2)^2(x^2+1)+(y-1)^2(y^2+1)
\end{eqnarray*}
%\end{example}
This example is a global unconstrained minimization problem. We take its gradient ideal and hereafter we consider the problem 
of minimizing the aforementioned function over its gradient ideal.
%\begin{example}
  \label{runex}
\begin{eqnarray*}
\inf_{\xx \in \R^2} \ f(x,y)=(x-1)^2(x-2)^2(x^2+1)+(y-1)^2(y^2+1)\\
 s.t. \ 6x^5-30x^4+56x^3-54x^2+34x-12=0 \nonumber \\
   4y^3-6y^2+4y-2=0
\end{eqnarray*}
The minimizer points are $(1,1)$ and $(2,1)$. The minimum is $f^*=0$.
\end{example}

\section{Convex relaxations}

In this section, we describe the finite dimensional convex 
optimization problems that we consider to solve the polynomial
optimization problem \eqref{problem1}. We recall the well-known 
full moment matrix relaxation and then we explain the border basis relaxation 
that we use. At the end of the section we compute the border basis for our running example.

But first, we introduce the notation we are going to use.
%\noindent{}\textbf{Notation.}
Let $\R[\xx]$ be the set of the polynomials in the variables
$\xx=(x_1,\ldots$, $x_n)$, with real coefficients in $\R$.
%in the field $\R$. %Hereafter, we will choose
%\footnote{For notational simplicity, we will consider only these two fields in this paper, but $\R$ and $\C$ can
%be replaced respectively by any real closed field and any field containing its algebraic closure.}
 %$\K=\R$. %or $\C$.
%Let $\overline{\K}$ denote the algebraic closure of ${\K}$.
For $\alpha \in \N^n$, $\xx^{\alpha}= x_1^{\alpha_1} \cdots x_n^{\alpha_n}$ is the monomial with exponent $\alpha$
and degree $|\alpha|=\sum_i\alpha_i$.  The set of all monomials in $\xx$ is
denoted $\Mon = \Mon(\xx)$. 
% We say that $\xx^{\alpha} \le \xx^{\beta}$ if
% $\xx^{\alpha}$ divides $\xx^{\beta}$, i.e., if $\alpha\le \beta$ coordinate-wise.
For a polynomial $f=\sum_\alpha f_\alpha \xx^\alpha$, its support is 
$supp(f):=\{\xx^\alpha\mid f_\alpha\ne 0\}$, the set of monomials occurring with a nonzero coefficient in $f$.

 For $t\in \N$ and $F\subseteq \R[\xx]$, we introduce the following sets:
 $F_{t}$ is the set of elements of $F$ of degree $\le t$; 
 $\Span{F} = \big\{ \sum_{f\in F} \lambda_{f}\, f\ |\ f\in F,
 \lambda_f\in \R\big\}$ is the linear span of $F$;
if $F$ is a vector space, ${F}^{*}$ is the dual space of linear forms from ${F}$ to $\R$;
 $(F) = \big\{ \sum_{f\in F} p_f\, f \ | \ p_f \in \R[\xx], f \in F \big\}$ is the ideal in $\R[\xx]$ generated by $F$;
 $\SpanD{F}{t}$ is the vector space spanned by $\{\xx^\alpha  f\mid f\in F_t, |\alpha|\le t-\deg(f)\}$; 
$F\cdot F:=\{p\,q\mid p,q\in F\}$;
$\Sigma^{2} (F) = \{ \sum_{i=1}^{s} f_{i}^{2}\mid f_{i} \in F\}$ is
the set of finite sums of squares of elements of $F$;
for $F=\{f_{1},\ldots,f_{m}\}\subset \R[\xx]$, 
$\prod (F) = \{ \prod_{i=1}^{m} f_{i}^{\epsilon_{i}}\mid
\epsilon_{i} \in \{0,1\} \}$.

%Hereafter, we will consider $\K=\R$.

\subsection{Hierarchies of relaxation problems}

\begin{defn}
 Given a finite dimensional vector space $E \subset \R[\xx]$ and a set of constraints $G
 \subset \R[\xx]$,  we define  the quadratic module of $G$ on $E$ as
$$
\begin{array}{l}
\Qc_{E,G} = \{ \sum_{g\in G^{^{0}}} g\, h+ \sum_{g'\in G^{+}} g'\, h'\\
 \ \ \ \ \mid h \in E, g\, h \in \Span{E\cdot E}, h' \in \Sigma^{2} (E),
 g'h'  \in \Span{E\cdot E} \}.
\end{array}
$$
\end{defn}
If $G^{\star} \subset \R[\xx]$ is the set of constraints such that
$G^{\star 0}=G^{0}$ and $G^{\star +}=\prod (G^{+})$,
the (truncated) quadratic module $\Qc_{E,G^{\star}}$ is called the
(truncated) preordering of $G$
and denoted $\Qc^{\star}_{E,G} $ or $\Pc_{E,G}$.

By construction,  $\Qc_{E,G} \subset \Span{E \cdot E}$ is a cone of
polynomials which are non-negative on the semi-algebraic set $S$.

We consider now its dual cone.
\begin{defn}
 Given a finite dimensional vector space $E \subset \R[\xx]$ which
 contains $1$ and a set of constraints $G
 \subset \R[\xx]$, we define 
$$
  \Lc_{E,G}:=\{\Lambda \in \Span {E \cdot E}^{*} \mid \Lambda(p) \geq
  0,\ \forall p \in \Qc_{E,G}, \Lambda (1)=1 \}.
$$
\end{defn}
The convex set associated to the preordering ${\Qc}^{\star}_{E,G}= {\Qc}_{E,G^{\star}}$ is
denoted ${\Lc}^{\star}_{E,G}$.

%The set $\Lc_{E,G}$ is the intersection of the closed convex cone of
%semi-definite positive quadratic forms on $E\times E$ with a linear
%space, thus it is a convex closed semi-algebraic set. 

By this definition, for any element $\Lambda \in \Lc_{E,G}$ and any $g\in
\Span{G^{0}}\cap E$, we have $\Lambda (g)=0$.

We introduce now truncated Hankel operators, which will play a central
role in the construction of the minimizer ideal of $f$ on $S$.
\begin{defn}
For a linear form $\Lambda \in \Span{E \cdot E}^{*}$, we define the map $H_{\Lambda}^E : E \rightarrow
 E^{\ast}$ by $H_{\Lambda}^E(p) (q) = \Lambda (p\,q)$ for $p,q\in E$.
 It is called the truncated Hankel operator of $\Lambda$ on the subspace $E$.
\end{defn}
Its matrix in the monomial and dual  bases of $E$ and $E^{*}$ is usually called the moment
matrix of $\Lambda$.
The kernel of this truncated Hankel operator %:
 % \begin{equation*}
 %  \ker H_{\Lambda}^E=\{p \in E \mid p \cdot \Lambda =0, i.e, \ \Lambda(pq)=0 \ \forall q \in E \}.
 % \end{equation*}
will be used to compute generators of the minimizer ideal, as we will see.
\begin{defn}\label{mintruncated} %\label{DefG}
 Let $E \subset \R[\xx]$ such that $ 1 \in E$ and a set of constraints
 $G \subset \R[\xx]$. We define the following extrema:
  \begin{itemize}
  %\item $f^*= \inf_{\xx \in S} f(\xx),$
  \item $f^{\mu}_{E,G}= \inf \ \{ \Lambda(f)$ s.t.  $\Lambda
    \in\Lc_{E,G} \},$
  \item $f^{sos}_{E,G}= \sup \ \{ \gamma \in \R$ s.t. $f-\gamma \in  \Qc_{E,G}\}.$  
 \end{itemize}
 By convention if the sets are empty, $f^{sos}_{E,G}=-\infty$ and $f^{\mu}_{E,G}=+\infty$.
\end{defn}
 If $E=\R[\xx]_t$ and $G^{0}=\SpanD{\gb^0}{2t}$, we also denote
$f^{\mu}_{E,G}$ by $f^{\mu}_{t,G}$ and $f^{sos}_{E,G}$ by $f^{sos}_{t,G}$.

We easily check that $f^{sos}_{E,G} \le f^{\mu}_{E,G}$, since if there exists $\gamma \in \R$ such that $f- \gamma = q \in \Qc_{E,G}$ 
 then $\forall \Lambda \in \mathcal{L}_{E,G}$, $\Lambda(f-\gamma)=
 \Lambda(f)- \gamma =  \Lambda(q) \ge 0$. 

If $\Sc (G)\subset S$, we
 also have $f^{\mu}_{E,G}\le f^{*}$ since for any $\mathbf{s} \in S$, the
 evaluation ${\mathbf{1}}_{\mathbf{s}}: p\in \R[\xx]\mapsto p (\mathbf{s})$ is 
 in $\Lc_{E,G}$.

Notice that if $E \subset E'$, $G \subset G'$ then $\Qc_{E,G} \subset
\Qc_{E',G'}$, $\mathcal{L}_{E',G'} \subset \mathcal{L}_{E,G}$,
$f^{\mu}_{E,G} \le f^{\mu}_{E',G'}$ {and} $f^{sos}_{E,G} \le f^{sos}_{E',G'}$.

\subsection{Full moment matrix relaxation hierarchy}
The relaxation hierarchies introduced in \citep{Las01} 
correspond to the case where $E=\R[\xx]_{t}$, $G^{0} =
\Span{\gb^{0}| 2t}$ and $G^{+}=\gb^{+}$. 

The quadratic module $\Qc_{\R[\xx]_{t},G}$ is denoted $\Qc_{t,\gb}$ and 
$\Lc_{\R[\xx]_{t},G}$ is denoted $\Lc_{t,\gb}$.
Hereafter, we will also call the
Lasserre hierarchy, the full moment matrix relaxation hierarchy. It corresponds
to the sequences 
$$ 
\cdots \subset \mathcal{L}_{t+1,\gb} \subset
\mathcal{L}_{t,\gb} \subset \cdots 
\ \mathrm{and}\  
\cdots \subset \Qc_{t,\gb} \subset
\Qc_{t+1,\gb} \subset \cdots 
$$
which yield the following increasing
sequences for $t \in \N$:
$$
 \cdots f^{\mu}_{t,\gb} \le f^{\mu}_{t+1,\gb} \le \cdots \le f^*
 \ \mathrm{and}\   
\cdots f^{sos}_{t,\gb} \le f^{sos}_{t+1,\gb} \le \cdots \le f^*.
$$
%\color{red}
The foundation of Lasserre's method is to show that these sequences
converge to $f^{*}$. This is proved under some conditions in
\citep{Las01}. It has also been shown that the limit can even be reached in a finite
number of steps in some cases, see e.g. \citep{LLR08b,NDS,Marshall,Ha-Pham:10,Nie11,Cert}.
In this case, the relaxation is said to be {\em exact}.
%\color{black}
\subsection{Border basis relaxation hierarchy}

In the following we are going to use another type of
relaxation hierarchy, which involves border basis computation.
Its aim is to reduce the size of the convex optimization problems 
solved at each level of the relaxation hierarchy. 
As we will see in Section \ref{sec:experiment}, the impact on the
performance of the relaxation approach is significant.
We briefly recall the properties of border basis that we need and
describe how they are used in the construction of this relaxation hierarchy.

Given a vector space $E \subseteq \R[\xx]$, its prolongation $ E^+ : =
E + x_1 E + \ldots + x_n E$ is again a vector space.

The vector space $E$ is said to be \textit{connected to 1} if $1\in E$
and there exists a finite increasing sequence of vector spaces $E_{0}\subset
E_{1} \subset \cdots \subset E$ such that $E_{0}=\Span{1}$, $E_{i+1} \subset E_{i}^{+}$.
For a monomial set $B\subseteq \Mon$, $B^{+}= B\cup x_{1} B \cup
\cdots \cup x_{n} B$ and $\partial B = B^{+} \setminus B$.
We easily check that $\Span{B}^{+}= \Span{B^{+}}$ and $\Span{B}$ is
connected to $1$ iff $1\in B$ and for every monomial $m\neq 1$ in $B$, $m= x_{i_{0}} m'$ for
some $i_{0}\in [1,n]$ and some monomial $m'\in B$. 
In this case, we will say that the monomial set $B$ is connected to $1$.

\begin{defn} \label{def:borderbasis}
Let $B \subset \Mon$ be connected to $1$. 
A family $F\subset \R[\xx]=R$ is a border basis for $B$ in degree $t\in \N$,
if $\forall f,f' \in F_{t}$,
\begin{itemize}
\item $supp(f)\subseteq B^+ \cap R_{t}$,
\item $f$ has exactly {\bf one} monomial in $\partial B$, denoted $\gamma(f)$ and called the leading monomial of f.
\item $\gamma(f)=\gamma(f')$ implies $f=f'$,
\item $\forall m\in \partial B \cap R_{t}$, 
  $\exists\, f \in F$ s.t. $\gamma (f) = m$,
\item $R_{t} = \Span{B}_{t} \oplus \Span {F | t}$.
\end{itemize}
A border basis $F$ for $B$ in all degrees $t$ is called a
border basis for $B$.
$F$ is \emph{graded} if moreover $\deg (\gamma(f)) = \deg (f) \ \forall  f\in F$.
\end{defn}
There are efficient algorithms to check that a given family $F$
is a border basis for $B$ in degree $t$ and to construct such family
from a set of polynomials. We refer to \citep{m-99-nf,Mourrain2005,BMPhT08,BMPhT12} for more details. We will use these
tools as ``black boxes'' in the following.  

For a border basis $F$ for $B$ in degree $t$, 
we denote by $\pi_{F,B}$ the projection of $R_{t}$ on $\Span{B_{t}}$
along $\Span {F | t}$. We easily check that 
\begin{itemize}
 \item $\forall m \in B_{t}$, $\pi_{B,F}(m)=m$,
 \item $\forall m \in \partial B \cap R_{t}$,
 $\pi_{B,F}(m)=m-f$, where $f$ is the (unique) polynomial in $F$ for which
 $\gamma(f)=m$, assuming the polynomials $f\in F$ are
normalized so that the coefficient of $\gamma(f)$ is $1$.
\end{itemize}

If $F$ is a graded border basis in degree $t$, one easily verifies that
$\deg(\piFB(m))\le \deg(m)$ for $m\in\Mon_t$. 

\textbf{Border basis hierarchy.} 
The sequence of relaxation problems that we will use hereafter is
defined as follows.
For each $t\in \N$, we construct the graded border basis $F_{2t}$ of $\gb^{0}$  in
degree $2\,t$. Let $B$ be the set of monomials (connected to
$1$) for which $F$ is a border basis in degree $2t$. We define
%\begin{itemize}
% \item 
$E_{t} := \Span{B_{t}}$,
% \item 
$G_{t}$ is the set of constraints such that $G_{t}^{0}= \{m
   -\pi_{B_t,F_{2t}} (m)$, $m \in B_{t} \cdot B_{t}\}$ and $G_{t}^{+}=\pi_{B_t,F_{2t}}(\gb^{+})$,
%\end{itemize}
and consider the relaxation sequence 
\begin{equation}\label{relaxation:bb}
\Qc_{E_{t}, G_{t}} \subset \Span{B_{t}\cdot B_{t}} \ \mathrm{and}\ 
\Lc_{E_{t}, G_{t}} \subset \Span{B_{t}\cdot B_{t}}^{*}  \
\end{equation}
{for} $t\in \N$. 
Since the subsets $B_{t}$ are not necessarily nested, these convex sets are
not necessarily included in each other. However,
by construction of the graded border basis of $\gb$, we have the following inclusions 
$$
\cdots \subset \Span{F_{2t}|2t} \subset  \Span{F_{2t+2}|2t +2} \subset
\cdots (\gb^{0}),
$$
and we can relate the border basis relaxation sequences with the
corresponding full moment matrix relaxation hierarchy, 
using the following proposition:
\begin{prop}\label{lemextend}
Let $t\in \N$, $B\subset \R[\xx]_{2t}$ be a monomial set connected to
$1$, 
$F\subset \R[\xx]$ be a border basis for $B$ in degree $2t$,
$E := \Span{B_{t}}$,
$E' := \R[\xx]_{t}$, 
$G,G'$ be sets of constraints such that 
$G^{0}= \{m -\pi_{B,F} (m)$, $m \in B_{t} \cdot B_{t}\}$, 
$G'^{0} = \Span {F| 2\,t}$,
$G^{+} = G'^{+}$.
Then for all
$\Lambda \in \Lc_{E,G}$, there exists a unique
$\Lambda' \in \Lc_{E',G'}$ which extends $\Lambda$.
Moreover, $\Lambda'$ satisfies
$\rank\, H_{{\Lambda'}}^{E'} = \rank\, H_{{\Lambda}}^{E}$ and
$\ker H_{{\Lambda'}}^{E'} = \ker H_{{\Lambda}}^{E} + \SpanD{F}{t}$.
\end{prop}
\begin{Proof}
As $F\subset \R[\xx]$ is a border basis for $B$ in degree $2t$, we have
 $\R[\xx]_{2\,t}=\Span{B}_{2\,t} \oplus \SpanD{F}{2\,t}$.
 As $\Span{B_{t}\cdot B_{t}} \subset \Span{B}_{2\,t} \oplus  \Span{G^{0}}$,
 $\Span{G^{0}} \subset \Span{G'^{0}} = \SpanD{F}{2t}$ and 
$\R[\xx]_{2\,t}= \Span{B}_{2\,t} \oplus \SpanD{F}{2t}$, we deduce
that for all $\Lambda\in \Lc_{E,G}$,
there exists a unique ${\Lambda'}\in  \R[\xx]_{2\,t}^{*}$
s.t. $\Lambda'_{|\Span{B}_{2t}}= \Lambda$ and $\Lambda' (\SpanD{F}{2t})=0$.

Let us first prove that $\Lambda' \in \Lc_{E',G'}= \Lc_{t,G'}$.
As any element $q'$ of $\Qc_{E',G'}$ can be decomposed as a sum of 
an element $q$ of $\Qc_{E,G}$ 
and an element $p\in \Span{F|2t}$, we have 
$\Lambda' (q') = \Lambda' (q) +\Lambda' (p) = \Lambda (q) \ge 0$.
This shows that $\Lambda' \in \Lc_{E',G'}$.

Let us prove now that  $\ker H_{{\Lambda'}}^{E'} = \ker H_{{\Lambda}}^{E}
+ \SpanD{F}{t}$ where $E := \Span{B_{t}}$, $E' := \R[\xx]_{t}$.
As $E\cdot \SpanD{F}{t} \subset \SpanD{F}{2t}=G'^{0}$, we have 
${\Lambda'} ( E\cdot  \SpanD{F}{t} ) =0$ so that 
\begin{equation}\label{eq:incl1}
\SpanD{F}{t} \subset \ker H_{{\Lambda'}}^{E'}.
\end{equation}

For any element $b\in \ker H_{\Lambda}^{E}$ we have
$\forall b'\in E,$ $\Lambda (b\,b')={\Lambda'} (b\,b') = 0$. As
${\Lambda'} (E\cdot \SpanD{F}{t})=0$ and $E'= E \oplus
\SpanD{F}{t}$, for any element $e\in E$, ${\Lambda'}(b\,e) = 0$. This proves that 
\begin{equation}\label{eq:incl2}
\ker H_{{\Lambda}}^{E}\subset \ker H_{{\Lambda'}}^{E'}.
\end{equation}
 
Conversely as $E'= E \oplus \SpanD{F}{t}$, any element
of $E'$ can be reduced modulo $\SpanD{F}{t}$ to an
element of $E$, which shows that 
\begin{equation}\label{eq:incl3}
\ker H_{{\Lambda'}}^{E'} \subset \ker H_{{\Lambda}}^{E} + \SpanD{F}{t}.
\end{equation}
From the inclusions \eqref{eq:incl1}, \eqref{eq:incl2} and \eqref{eq:incl3}, we deduce that $\ker
H_{{\Lambda'}}^{E'}$ = $\ker H_{{\Lambda}}^{E} + \SpanD{F}{t}$
and that $\rank\, H_{{\Lambda'}}^{E'} = \rank\, H_{{\Lambda}}^{E}$.
\end{Proof}

We deduce from this proposition that $f_{E_{t},G_{t}}^{\mu}=f_{t,\Span{F_{2t}| 2t}}^{\mu}$.
The sequence of convex sets $\Lc_{E_{t}, G_{t}}$ can be seen 
as the projections of nested convex sets 
$$
\cdots
\supset
\Lc_{t,\gb} \supset 
\Lc_{t+1,\gb} \supset 
%\supset 
\cdots
%\end{array}
$$
so that we have
$\cdots \le f^{\mu}_{E_{t},G_{t}} \le f^{\mu}_{E_{t+1},G_{t+1}} \le \cdots \le f^*$.
We check that similar properties hold for  $\Qc_{E_{t}, G_{t}}$, $\Qc_{t,\gb}$ and
$f_{E_{t},G_{t}}^{sos}=f_{t,\gb}^{sos}$, taking the
quotient modulo $\Span{F_{2t}| 2t}$.

\vspace{1cm}
Now we compute the border basis for our running example \ref{runex} and the monomials that
we can reduce by using this border basis.
\begin{example} 
We take the set of constraint $g^0=\{6x^5-30x^4+56x^3-54x^2+34x-12,4y^3-6y^2+4y-2\}$
and $t=3$.
 \begin{itemize}

  \item The border basis is $F_3=\{x^5-5x^4+9.333x^3-9x^2+5.66x-2,y^3-1.5y^2+y-0.5\}$
  \item The monomial basis in degree $\leq 3$ is:
  \begin{equation*}
   B_3=\{1,x,y,x^2,xy,y^2,x^3,x^2y,xy^2\}
  \end{equation*}
   The monomial $y^3$ is the leading term of an element of $F_{3}$.
  \item The border basis SDP relaxation is constructed from the
    reduction of the monomials in $B_{3}\cdot B_{3}$. The following
    monomials are reduced by the border basis
  \begin{equation*}
   \{y^3,xy^3,y^4,x^5,x^2y^3,xy^4,x^6,x^5y,x^3y^3,x^2y^4\}.
  \end{equation*}
  This yields the following constraints:
  \begin{align*}
& y^3\equiv 0.5-\,y+1.5\,y^2 \\  
& xy^3\equiv 0.5\,x-\,xy+1.5\,xy^2 \\   
& y^4\equiv 0.75-\,y+1.25\,y^2 \\   
& x^2y^3\equiv 0.5\,x^2-\,x^2y+1.5\,x^2y^2\\  
& xy^4\equiv 0.75\,x-\,xy+1.25\,xy^2 \\   
& x^5\equiv 2-5.666\,x+9\,x^2-9.333\,x^3+5\,x^4\\   
& x^6 \equiv 10-26.333\,x+39.333\,x^2-37.666\,x^3+15.666\,x^4 \\   
& x^5y\equiv 2\,y-5.666\,xy+9\,x^2y-9.333\,x^3y+5\,x^4y \\
& x^3y^3\equiv 0.5\,x^3-\,x^3y+1.5\,x^3y^2 \\   
& x^2y^4\equiv 0.75\,x^2-\,x^2y+1.25\,x^2y^2 \\ 
\end{align*}

 \end{itemize}

\end{example}

\section{Optimal linear form}\label{sec:3}
In this section we introduce the notion of {\em optimal  linear form for $f$},
involved in the computation of $I_{min}$ (also called 
generic linear form when $f=0$ in \citep{LLR08b, lasserre:hal-00651759}). In order to find this optimal linear form
we solve a Semi-Definite Programming (SDP) problem,
which involves truncated Hankel matrices associated with the monomial
basis and the reduction of their products 
by the border basis, as described in the previous section. This allows us to reduce the size of the matrix and the number of parameters. At the end of this section we compute the optimal linear form for our running example.

\begin{defn}\label{generic}
$\Lambda^* \in \mathcal{L}_{E,G}$ is optimal  for $f$ if $$
\rank \ H_{\Lambda^*}^{E} = \max_{\Lambda \in
  \mathcal{L}_{E,G},\Lambda(f)=f^{\mu}_{E,G}} \rank \ H_{\Lambda}^{E}.$$
\end{defn}
The next result shows that only elements in $I_{min}$ are involved in the kernel of a truncated
Hankel operator associated with an optimal  linear form for $f$.
\begin{thm}\label{kermin}
 Let $E \subset \R[\xx]$ such that $ 1 \in E$ and $f \in \langle E
 \cdot E \rangle$ and let $G \subset \R[\xx]$ be a set of constraints with
 $V_{min}\subset \Sc (G)$. If $\Lambda^* \in \Lc_{E,G}$ 
is optimal  for $f$ and such that $\Lambda^{*} (f) = f^{*}$, then $\ker
H_{\Lambda^*}^{E} \subset I_{min}$.  
\end{thm}
\begin{Proof}
The proof is similar e.g. to \citep{lasserre:hal-00651759}[Theorem 4.9].
\end{Proof}\\
Let us describe how optimal linear forms are computed by solving 
convex optimization problems:\\
\begin{algorithm2e}[H]\caption{\textsc{Optimal Linear
      Form}}\label{optl}
\noindent{}\textbf{Input:} $f\in \R[\xx]$, $B_t=(\xx^{\alpha})_{\alpha \in A}$ a monomial set of degree $\le t$ containing $1$ 
with $f= \sum_{\alpha \in A+ A} f_{\alpha} \xx^{\alpha} \in
\Span{B_t\cdot B_t}$, $G\subset \R[\xx]$.\\
\noindent{}\textbf{Output:} 
the infimum $f^{\mu}_{t,G}$ of $\sum_{\alpha\in A+ A} \lambda_{\alpha} f_{\alpha}$ subject  to:
     \begin{itemize}
        \item[--] $H_{\Lambda^*}^{B_t}= (h_{\alpha,\beta})_{\alpha,\beta\in A} \succcurlyeq 0$,
        \item[--] $H_{\Lambda^*}^{B_t}$ satisfies the Hankel constraints \\
               $h_{0,0}=1$, and $h_{\alpha,\beta}=h_{\alpha',\beta'}$
               if $\alpha+\beta=\alpha'+\beta'$,
       \item[--]  $\Lambda^{*} (g^{0})= \sum_{\alpha \in  A+
            A}\, g^{0}_{\alpha} \lambda_{\alpha}=0$ for all
          $g^{0}=\sum_{\alpha \in  A+ A}\, g^{0}_{\alpha} \xx^{\alpha}
          \in G^{0}\cap \Span{B_t\cdot B_t}$.
       \item[--] $H_{g^{+}\cdot\Lambda^*}^{B_{t-w}} \succcurlyeq 0$ for all $g^{+} \in G^{+}$  where $w=\lceil \frac{deg(g^{+})}{2} \rceil$.
     \end{itemize}
and $\Lambda^{*}\in \Span{B_t \cdot B_t}^{*}$
represented by the vector $[\lambda_{\alpha}]_{\alpha \in A + A}$.
\end{algorithm2e}
This optimization algorithm involves a Semidefinite programming problem,
corresponding to the optimization of a linear functional on the
intersection of a linear subspace with the convex set of positive semidefinite matrices.
It is a convex optimization problem, which can be solved efficiently
%(in polynomial time up to some precision \citep{XXX}) 
by SDP solvers.
If an interior point method is used, the solution $\Lambda^{*}$ is in the interior of
a face on which the infimum $\Lambda^{*} (f)$ is reached so that $\Lambda^{*}$ is optimal  for $f$.
This is the case for tools such as \textsc{csdp}, \textsc{sdpa}, \textsc{sdpa-gmp}, and \textsc{mosek}
that we will use in the experiments.

%\vspace{1cm}
%Hereafter, we solve the SDP problem associated to our running example \ref{runex}.
\begin{example}
  For the running example \ref{runex} and  the relaxation order $t=3$, we solve the following SDP problem:
\begin{align*} 
 &\inf \ \Lambda(f)=2.75-\Lambda(y)-4.333\Lambda(x)+0.25\Lambda(y^2)+2.333\Lambda(x^2)+0.333\Lambda(x^3)-0.333\Lambda(x^4) \\
 & with \ \Lambda \ s.t. \\
&\Lambda(y^3)=0.5-\Lambda(y)+1.5\Lambda(y^2) \\  
&\Lambda(xy^3)=0.5\Lambda(x)-\Lambda(xy)+1.5\Lambda(xy^2) \\   
&\Lambda(y^4)=0.75-\Lambda(y)+1.25\Lambda(y^2) \\   
&\Lambda(x^5)=2-5.666\Lambda(x)+9\Lambda(x^2)-9.333\Lambda(x^3)+5\Lambda(x^4)\\   
&\Lambda(x^2y^3)=0.5\Lambda(x^2)-\Lambda(x^2y)+1.5\Lambda(x^2y^2)\\  
&\Lambda(xy^4)=0.75\Lambda(x)-\Lambda(xy)+1.25\Lambda(xy^2) \\   
&\Lambda(x^6)=10-26.333\Lambda(x)+39.333\Lambda(x^2)-37.666\Lambda(x^3)+15.666\Lambda(x^4)  \\
& \Lambda(x^5y)=2\Lambda(y)-5.666\Lambda(xy)+9\Lambda(x^2y)-9.333\Lambda(x^3y)+5\Lambda(x^4y) \\      
&\Lambda(x^3y^3)=0.5\Lambda(x^3)-\Lambda(x^3y)+1.5\Lambda(x^3y^2) \\   
&\Lambda(x^2y^4)=0.75\Lambda(x^2)-\Lambda(x^2y)+1.25\Lambda(x^2y^2) \\
&\Lambda (1)=1 
\end{align*}
and 
\begin{align*}
H_{\Lambda}^{B_3} :=\left(
\begin{array}{cccccccccc}
1& a & b & c & d & e & f & g & h  \\
a & c & d & f& g & h & i & j & k \\
b & d & e & g & h & \Lambda(y^3) & j & k & \Lambda(xy^3) \\
c & f & g & i & j & k & \Lambda(x^5) & l & m \\
d & g & h & j & k & \Lambda(xy^3) & l & m & \Lambda(x^2y^3) \\
e & h & \Lambda(y^3) & k & \Lambda(xy^3) &\Lambda(y^4) & m &\Lambda(x^2y^3) & \Lambda(xy^4)\\
f & i & j & \Lambda(x^5)& l &m&\Lambda(x^6) & \Lambda(x^5y) &n \\
g & j & k & l & m & \Lambda(x^2y^3) & \Lambda(x^5y)& n &\Lambda(x^3y^3) \\
h & k & \Lambda(xy^3) & m & \Lambda(x^2y^3) &\Lambda(xy^4) &n & \Lambda(x^3y^3) & \Lambda(x^2y^4)
\end{array}
\right) \succcurlyeq 0
\end{align*}
where  $a=\Lambda (x)$, $b=\Lambda (y)$, $c=\Lambda (x^2)$, $d=\Lambda (xy)$, $e=\Lambda (y^2)$, $f=\Lambda (x^3)$,\\
$g=\Lambda (x^2y)$, $h=\Lambda (xy^2)$, $i=\Lambda (x^4)$,$j=\Lambda (x^3y)$, $k=\Lambda (x^2y^2)$, $l=\Lambda (x^4y)$,\\
$m=\Lambda (x^3y^2)$, $n=\Lambda (x^4y^2)$ and 
$\Lambda(y^3)=0.5-b+1.5e, \Lambda(y^4)=0.75-b+1.25e,\\
\Lambda(x^2y^3)=0.5c-g+1.5k,\Lambda(xy^3)=0.5a-d+1.5h,\\
\Lambda(x^5)=2-5.666a+9c-9.333f+5i,\Lambda(x^5y)=2b-5.666d+9g-9.333j+5l,\\
\Lambda(x^6)=10-26.333a+39.333c-37.666f+15.666i,\Lambda(xy^4)=0.75a-d+1.25h,\\
\Lambda(x^3y^3)=0.5f-j+1.5m,\Lambda(x^2y^4)=0.75c-g+1.25k.$
\vspace{0.5cm}

A \textbf{solution} is:
$\Lambda^*(1)=1,\Lambda^*(x)=1.5,\Lambda^*(y)=1,\Lambda^*(x^2)=2.5,\Lambda^*(xy)=1.5,\Lambda^*(y^2)=1,\Lambda^*(xy^2)=1.5,\Lambda^*(x^2y)=2.5,\Lambda^*(x^3)=4.5,
\Lambda^*(x^2y^2)=2.5,\Lambda^*(x^3y)=4.5,\Lambda^*(x^4)=8.5,\Lambda^* (x^4y)=4.5,\Lambda^* (x^3y^2)=8.5, \Lambda^* (x^4y^2)=8.5$.\\
The minimum is $\Lambda^* (f)=0$.
\end{example}

%\section{Convergence certification}\label{sec-0dim}
\section{Decomposition}\label{sec:4}

To be able to compute the minimizer points from an optimal linear
form, we need to detect when the infimum is reached.
In this section, we describe new criterion to check when the kernel of
a truncated Hankel operator associated to an optimal linear form for
$f$ yields the generators of the minimizer ideal. It involves the 
flat extension theorem of \citep{MoLa2008} and applies to polynomial optimization
problems where the minimizer ideal $I_{min}$ is zero-dimensional.
At the end of this section we verify the flat extension property in our running example.

\subsection{Flat extension criterion}

\begin{defn} \label{defflatextension}
Given vector subspaces $E_0 \subset E \subset \R[\xx]$ and $\Lambda \in { \Span{E \cdot E}}^*$, $H_{\Lambda}^E$ is
said to be a \textit{flat extension} of its restriction $H^{E_0}_{\Lambda}$ if $\rank\, H^E_{\Lambda} = \rank\, H^{E_0}_{\Lambda}$.
\end{defn}

We recall here a result from \citep{MoLa2008}, which gives a rank condition
for the existence of a flat extension of a truncated Hankel operator\footnote{In \citep{MoLa2008}, it is stated with a vector space spanned by a monomial set connected
to $1$, but its extension to vector spaces connected to $1$ is straightforward.}.
\begin{thm}\label{theoflatextension}
Let $V\subset E \subset \R[\xx]$ be vector spaces connected to 1 with 
$V^{+} \subset E$ and let  $\Lambda\in \Span{E \cdot E}^{*}$. 
Assume that $\rank\, H^{E}_{\Lambda} = \rank\, H^V_{\Lambda} = \dim V$. 
Then there exists a (unique) linear form $\tilde{\Lambda} \in \R[\xx]^*$ which extends
 $\Lambda$, i.e., $\tilde{\Lambda} (p) = \Lambda (p)$ for all $p \in
 \Span{E \cdot E}$, 
satisfying $\rank\, H_{\tilde{\Lambda}} = \rank\, H^{E}_{\Lambda}$. Moreover, we have $\ker H_{\tilde{\Lambda}}= (\ker H^{E}_{\Lambda})$.
\end{thm}
In other words, the condition $\rank\, H^{E}_{\Lambda} = \rank\,
H^V_{\Lambda} = \dim V$ 
implies that the truncated Hankel operator $H^{E}_{\Lambda}$ has
a (unique) flat extension to a (full) Hankel operator
$H_{\tilde{\Lambda}}$ defined on $\R[\xx]$.

\begin{thm}\label{caszero}
Let $V\subset E \subset \R[\xx]$ be finite dimensional vector spaces connected to $1$
with $V^{+} \subset E$,
$G^{0}\cdot V\subset \Span{E\cdot E}$,
$G^{+}\cdot V \cdot V \subset \Span{E\cdot E}$.

Let $\Lambda\in \mathcal{L}_{E,G}$ such that
$\rank\, H_{\Lambda}^{E} = \rank\, H_{\Lambda}^{V}= \dim V $. 
Then there exists a linear form $\tilde{\Lambda} \in
   \R[\xx]^{*}$ which extends $\Lambda$ and is supported on points
   of $\Sc (G)$ with positive weights:
$$ 
\tilde{\Lambda} = \sum_{i=1}^{r} \omega_{i} \ev_{\xi_{i}}\ \mathrm{with}\
\omega_{i} >0, 
\xi_{i}\in \Sc (G).
$$
Moreover, $(\ker H_{\Lambda}^{E}) = \Ic(\xi_{1}, \ldots, \xi_{r})$.
\end{thm}

\begin{Proof}
As $\rank\, H_{\Lambda}^{E} = \rank\, H_{\Lambda}^{V}= \dim V$,
Theorem \ref{theoflatextension} implies that there exists a (unique)
linear function $\tilde{\Lambda} \in \R[\xx]^*$ which extends
$\Lambda$.
As $\rank\, H_{\tilde{\Lambda}} = \rank\,
 H^{V}_{\Lambda}=|V|$ and $\ker H_{\tilde{\Lambda}}= (\ker
 H^{E}_{\Lambda})$, any polynomial $p\in \R[\xx]$ can be reduced modulo
 $\ker H_{\tilde{\Lambda}}$ to a polynomial $b\in V$ so
 that $p-b\in \ker H_{\tilde{\Lambda}}$.
Then $\tilde{\Lambda} (p^{2})= \tilde{\Lambda} (b^{2}) =
\Lambda (b^{2}) \geq 0$ since $\Lambda\in \mathcal{L}_{E,G}$.
%This implies that $\tilde{\Lambda} \succcurlyeq 0$. 
By Theorem 3.14 of \citep{lasserre:hal-00651759},
$\tilde{\Lambda}$ has a decomposition of the form
$\tilde{\Lambda}=\sum_{i=1}^r \omega_i\ev_{\xi_i}$  
with $\omega_i > 0$ and $\xi_i \in \R^{n}$.

By Lemma 3.5 of \citep{lasserre:hal-00651759}, $V$ is isomorphic to $\R[\xx]/\Ic(\xi_{1}, \ldots,
\xi_{r})$ and there exist (interpolation) polynomials $b_{1},
\ldots, b_{r} \in V$ satisfying $b_{i} (\xi_{j}) =1$ if $i=j$
and $b_{i} (\xi_{j}) =0$ otherwise. We deduce that for $i=1, \ldots,
r$ and for all elements $g\in G^{0}$, 
$$
\Lambda (b_{i} g) = 0 = \tilde{\Lambda} (b_{i} g) = \omega_i g
(\xi_{i}).$$ 
As $\omega_i > 0$ then $g(\xi_{i})=0$.
Similarly, for all $h\in G^{+}$, 
$$\Lambda (b_{i}^{2} h) = 
\tilde{\Lambda} (b_{i}^{2} h) = \omega_i h (\xi_{i}) \geq 0
$$ and
$h (\xi_{i})\geq 0$, hence
$\xi_{i} \in \Sc (G)$.

By Theorem 3.14 of \citep{lasserre:hal-00651759} and Theorem
\ref{theoflatextension}, we also have
$\ker H_{\tilde{\Lambda}}
 = \Ic(\xi_{1}, \ldots, \xi_{r}) = (\ker H^{E}_{\Lambda}).$
\end{Proof}

This theorem applied to an optimal  linear form $\Lambda^{*}$ for $f$
gives a convergence certificate to check when the infimum $f^{*}$ is
reached and when a generating family of the minimizer ideal is
obtained.  It generalizes the flat truncation certificate given in
\citep{Nie2012}. As we will see in the experiments, it allows to
detect more efficiently when the infimum is reached.
Notice that if the test is satisfied, necessarily $I_{min}$ is zero-dimensional.

\subsection{Flat extension algorithm}

In this section, we describe a new algorithm to check the flat
extension property for a linear form for which some moments are known.

Let $E$ be a finite dimensional subspace of $\R[\xx]$ connected to $1$ and let
$\Lambda^{*}$ be a linear form defined on $\Span{E\cdot E}$ given 
by its ``moments'' $\Lambda^{*} (e_{i}):= \Lambda^{*}_{i}$, where 
$e_{1}, \ldots, e_{s}$ is a basis of $\Span{E\cdot E}$ (for instance a
monomial basis).
In the context of global polynomial optimization that we consider here, this linear form is
an optimal  linear form for $f$ (see Section \ref{sec:3})
computed by SDP.

We define the linear functional $\Lambda^{*}$ from its moments as
$\Lambda^{*} : p = \sum_{i=1}^{s} p_{i} e_{i}  \in \Span{E\cdot E}
\mapsto \sum_{i=1}^{s} p_{i} \Lambda_{i}$ and the corresponding inner product:
\begin{eqnarray}\label{eq:innerprod}
E \times E& \rightarrow & \R \nonumber\\
(p,q) & \mapsto & \Span {p,q}_{*}:= \Lambda^{*} (p\, q)
\end{eqnarray}

To check the flat extension property, we are going to inductively define
 vector spaces $V_i$ \ as follows. Start with $V_0 = \langle 1
\rangle$. Suppose $V_{i}$ is known and compute a vector space $L_i$ of
maximal dimension in $V_i^+$ such that $L_i$ is orthogonal to $V_i$: \ $\langle L_i, V_i \rangle_{*} =
  0$ and $L_i \cap \ker H_{\Lambda^{*}}^{V_i^+} = \{ 0 \}$.
Then we define $V_{i + 1} = V_i + L_i$. 

Suppose that $b_1, \ldots, b_{r_i}$ is an orthogonal basis of $V_i$:
$\langle b_i, b_j \rangle_{*} = 0$ if $i \neq j$ and $\langle b_i, b_i
\rangle_{*} \neq 0$ . Then $L_i$ can be constructed as follows: Compute
the vectors \
\[ b_{i, j} = x_j b_i - \sum_{k = 1}^{r_i} \frac{\langle x_j \nobracket b_i,
   b_k \rangle_{*} \nobracket}{\langle b_k, b_k \rangle_{*}} b_k, \]
generating $V_i^{\perp}$ \ in $V_i^+$ and extract a maximal orthogonal family
\ $b_{r_i + 1}, \ldots, b_{r_i + s}$ for the inner product $\langle ., .
\rangle_{*}$, that form a basis of $L_i$. This can be done for instance
by computing a QR decomposition of the matrix
$ [ \langle b_{i, j}, b_{i', j'} \rangle_{*}]_{1 \leqslant i, i'
   \leqslant r_i, 1 \leqslant j, j' \leqslant n}$.
The process can be repeated until either
\begin{itemize}
 
\item $V_{i}^{+}\not\subset E$ and the algorithm will stop and return
   \texttt{failed},
 \item or $L_{i}=\{0\}$ and $V_{i}^{+} = V_{i} \oplus \ker
 H_{\Lambda^{*}}^{V_i^+}$. In this case, the algorithm stops with \texttt{success}.
\end{itemize}
Here is the complete description of the algorithm:\\
\begin{algorithm2e}[H]\caption{\textsc{Decomposition}\label{algo:flatcrit}}
{\tmstrong{Input:}} a vector space $E$ connected to $1$ and a linear
form $\Lambda^{*}\in \Span{E\cdot E}^{*}$.\\
{\tmstrong{Output:}} \texttt{failed} or \texttt{success} with 
\begin{itemize}
 \item a basis $B = \{ b_1, \ldots, b_r \}\subset \R[\xx]$;
 \item the relations $x_k b_j - \sum_{i = 1}^{r}
    \frac{\langle x_k \nobracket b_j, b_i \rangle_{*}
      \nobracket}{\langle  b_i, b_i \rangle_{*}} b_i$, $j=1\ldots r$
    $k=1\ldots n$.
\end{itemize}
{\tmstrong{Begin}}
    \begin{enumerate}
      \item Take $B:= \{1\}$; \ $s:=1$; $r := 1$; \
      
      \item While $s > 0$ and $B^{+}\subset E$ do
      \begin{enumerate}
        \item compute $b_{j,k} : = x_k b_j - \sum_{i = 1}^{r}
        \frac{\langle x_k \nobracket b_j, b_i \rangle_{*}
        \nobracket}{\langle b_i, b_i \rangle_{*}} b_i$ for $j = 1,
        \ldots, r$, $k = 1, \ldots, n$;
        
        \item compute a maximal \ subset $B'=\{ b'_{1}, \ldots, b'_{s} \}$
        of $\Span{ b_{j,k} }$ of orthogonal vectors for the inner product
        $\langle ., . \rangle_{*}$ and let $B:= B\cup B'$,
      $s=\mid B' \mid$ and $r \plusassign s$;
      \end{enumerate}
     \item If $B^{+} \not \subset E$ then return \texttt{failed} \\
           else $(s=0)$ return \texttt{success}.
    \end{enumerate}  
{\tmstrong{End}}
%  {\tmstrong{Output:}} \texttt{failed} or \texttt{success} with 
% \begin{itemize}
%  \item a basis $B = \{ b_1, \ldots, b_r \}\subset \R[\xx]$;
%  \item the relations $x_k b_j - \sum_{i = 1}^{r}
%     \frac{\langle x_k \nobracket b_j, b_i \rangle_{*}
%       \nobracket}{\langle  b_i, b_i \rangle_{*}} b_i$, $j=1\ldots r$
%     $k=1\ldots n$.
%\end{itemize}
\end{algorithm2e}

Let us describe the computation performed on the moment matrix, during
the main loop of the algorithm.
At each step, the moment matrix of ${\Lambda^{*}}$ on $V_{i}^{+}$ is of the form
$$  
  H^{V_{i}^{+}}_{\Lambda^{*}}=\left[
    \begin{array}{c|c}
      H_{\Lambda^*}^{B_i,B_i} & H_{\Lambda^*}^{B_{i},\partial
  B_i}  \\ \hline
      H_{\Lambda^*}^{\partial B_i, B_i} & H_{\Lambda^*}^{\partial B_{i},\partial
  B_i} 
     \end{array}\right]
$$
where $\partial B_{i}$ is  a subset of $\{b_{i,j}\}$ such that $B_{i} \cup \partial B_{i}$ is a
basis of $\Span{B_{i}^{+}}$.
By construction, the matrix $H_{\Lambda^*}^{B_i,B_i}$ is diagonal
since $B_{i}$ is orthogonal for $\Span{\cdot, \cdot }_{*}$.
As the polynomials $b_{i, j}$ 
%$ : = x_j b_i - $ $ \sum_{k = 1}^{r} \frac{\langle x_j \nobracket b_i, b_k \rangle_{*}}{\langle b_k, b_k \rangle_{*}} b_k$
are orthogonal to $B_{i}$, we have $H_{\Lambda^*}^{B_{i},\partial
  B_i}= H_{\Lambda^*}^{\partial B_i, B_i}=0$. 
If $H_{\Lambda^*}^{\partial B_{i},\partial B_i}=0$ then the algorithm
stops with \texttt{success} and all the elements $b_{i,j}$ are in the
kernel of $H_{\Lambda^{*}}^{B_{i},B_{i}}$. Otherwise an orthogonal basis $b'_{1},
\ldots, b'_{s}$ is extracted. It can then be completed in a basis
of $\Span{b_{i,j}}$ so that the matrix $H_{\Lambda^*}^{\partial B_{i},\partial B_i}$
in this basis is diagonal with zero entries after the $(s+1)^{th}$ index.
In the next loop of the algorithm, the basis $B_{i+1}$ contains the 
maximal orthogonal family $b'_{1},\ldots, b'_{s}$ so that the matrix $ H_{\Lambda^*}^{B_{i+1},B_{i+1}}$
remains diagonal and invertible. 

\begin{prop}\label{proposition:algo3.1}
Let $\Lambda^{*} \in \Lc_{E,G}$ be optimal  for $f$.
If Algorithm \ref{algo:flatcrit} applied to $\Lambda^{*}$ and $E$ stops with \texttt{success}, then 
\begin{enumerate}
 \item there exists a linear form $\tilde{\Lambda} \in
   \R[\xx]^{*}$ which extends $\Lambda^{*}$ and is supported on
   points in $\Sc(G)$ with positive weights:
$$ 
\tilde{\Lambda} = \sum_{i=1}^{r} \omega_{i} \ev_{\xi_{i}}\ \mathrm{with}\
\omega_{i} >0, 
\xi_{i}\in \R^{n}.
$$
 \item $B=\{b_{1}, \ldots, b_{r}\}$ is a basis of
 $\Ac_{\tilde{\Lambda}}=\R[\xx]/I_{\tilde{\Lambda}}$ where
 $I_{\tilde{\Lambda}} = \ker H_{{\tilde{\Lambda}}}$,
 \item $x_k b_j - \sum_{i = 1}^{r_{}}
        \frac{\langle x_k \nobracket b_j, b_i \rangle_{*}}{\langle
          b_i, b_i \rangle_{*}} b_i$,
$j=1,\ldots,r$, $k=1,\ldots,n$ 
are generators of $I_{\tilde{\Lambda}}= \Ic(\xi_{1}, \ldots,
\xi_{r})$,
 \item $f_{E,G}^{\mu}= f^{*}$,
 \item $V_{min}= \{\xi_{1}, \ldots, \xi_{r}\}$.
\end{enumerate}
\end{prop}
\begin{Proof}
When the algorithm terminates with \texttt{success}, the set $B$ is
such that $\rank\,  H^{B^{+}}_{\Lambda^{*}} =  \rank\, H^{B}_{\Lambda^{*}}= |B|$.
By Theorem \ref{caszero}, there exists a linear form $\tilde{\Lambda} \in
   \R[\xx]^{*}$ extends $\Lambda^{*}$ and is supported on points in
   $\Sc (G)$ with positive weights:
$$ 
\tilde{\Lambda} = \sum_{i=1}^{r} \omega_{i} \ev_{\xi_{i}}\ \mathrm{with}\
\omega_{i} >0, 
\xi_{i}\in \Sc(G).
$$
This implies that $\Ac_{\tilde{\Lambda}}$ is of dimension $r$ and
that $I_{\tilde{\Lambda}} = \Ic(\xi_{1}, $ $\ldots, \xi_{r})$.
As $H^{B}_{\Lambda^{*}}$ is invertible, $B$ is a basis of
$\Ac_{\tilde{\Lambda}}$ which proves the second point.

Let $K$ be the set of polynomials $x_j b_i - \sum_{k = 1}^{r_{}}
        \frac{\langle x_j \nobracket b_i, b_k \rangle_{*}}{\langle b_k, b_k \rangle_{*}} b_k$.
If the algorithm terminates with \texttt{success}, we have $\ker
H^{B^{+}}_{\Lambda^{*}} =\Span{K}$ and by Theorem \ref{caszero}, we
deduce that $(K)= (\ker H^{B^{+}}_{\Lambda^{*}}) =
I_{\tilde{\Lambda}}$, which proves the third point.

As $\tilde{\Lambda} (1)= 1$, we have $\sum_{i=1}^{r} w_{i}=1$ and 
$$
\tilde{\Lambda} (f) = \sum_{i=1}^{r} \omega_{i} f (\xi_{i}) \ge f^{*}
$$
since  $\xi_{i} \in \Sc(G)$ and $f (\xi_{i}) \ge f^{*}$. 
The relation $f^{\mu}_{E,G} \le f^*$
implies that $f (\xi_{i}) = f^{*}$ for $i=1,\ldots, r$ and  the fourth
point is true: $f^{\mu}_{E,G} = f^*$.

As $f (\xi_{i}) = f^{*}$ for $i=1,\ldots, r$, we have $\{\xi_{1}, \ldots, \xi_{r}\} \subset V_{min}$.
By Theorem \ref{kermin}, the polynomials of $K$ are in $I_{min}$ so that
$V_{min} \subset \Vc (K) = \{\xi_{1}, \ldots, \xi_{r}\}$.
This shows that $V_{min} = \{\xi_{1}, \ldots, \xi_{r}\}$ and concludes
the proof of this proposition
\end{Proof}

\begin{example}\label{exampledal}
We apply  Algorithm \ref{algo:flatcrit} to our running example.
A solution of the SDP problem output by Algorithm \ref{optl} is:\\
$\Lambda^*(1)=1,\Lambda^*(x)=1.5,\Lambda^*(y)=1,\Lambda^*(x^2)=2.5,\Lambda^*(xy)=1.5,\Lambda^*(y^2)=1,\Lambda^*(xy^2)=1.5,\Lambda^*(x^2y)=2.5,\Lambda^*(x^3)=4.5,
\Lambda^*(x^2y^2)=2.5,\Lambda^*(x^3y)=4.5,\Lambda^*(x^4)=8.5,\Lambda^* (x^4y)=4.5,\Lambda^* (x^3y^2)=8.5, \Lambda^* (x^4y^2)=8.5$\\
We verify the flat extension criterion for $\R[\xx]_3$.
\begin{itemize}

   \item $B_{0}=\{1\}, \ \partial B_{0}=\{x,y\}, \ B_{0}^{+}=\{1,x,y\}$\\
   
    \[H_{\Lambda^*}^{B_{0}^{+}}=  \left( \begin{array}{ccc}
    1 & 1.5 & 1\\
    1.5 & 2.5 & 1.5\\
    1& 1.5& 1\end{array} \right) \longrightarrow H_{\Lambda^*}^{\{1,x-1.5,y-1\}}=\left( \begin{array}{ccc}
    1 & 0 & 0\\
    0 & 0.25 & 0\\
    0 & 0& 0 \end{array} \right) \]\\
    \[\rank \ H_{\Lambda^*}^{B_{0}^{+}}=2,\]\\
     \[ \{y-1\} \in \ker H_{\Lambda^*}^{B_{0}^{+}}, \ \{x-1.5\} \perp B_{0} \ and \ \{x-1.5\} \notin \ker H_{\Lambda^*}^{B_{0}^{+}},\]\\
     \[L_{0}=\{x-1.5\}.\]
    
    \vspace{0.5cm}

   \item $B_{1}=B_{0}\cup L_{0}=\{1,x-1.5\}, \ \partial B_{1}=\{y,x^2-1.5x,xy-1.5y\},$\\
     $ B_{1}^{+}=\{1, x-1.5 ,y,x^2-1.5x,xy-1.5y\}$\\
     
     \[H_{\Lambda^*}^{B_{1}^{+}}=  \left( \begin{array}{ccccc}
    1 & 0 & 1&0.25&0\\
    0 & 0.25 & 0&0.375&0.25 \\
    1 &0& 1&0.25&0\\\
    0.25&0.375&0.25&0.625&0.375\\
    0&0.25&0&0.375&0.25\end{array} \right) \longrightarrow H_{\Lambda^*}^{\tilde{B}}=\left( \begin{array}{ccccc}
    1 & 0 &0&0&0\\
    0 & 0.25&0&0&0\\
    0&0&0&0&0\\
    0&0&0&0&0\\
    0&0&0&0&0\end{array} \right)  \] 
   
    where $\tilde{B}=\{1,x-1.5,y-1,x^2-3x+2,xy-1.5y-x+1.5\}$\\
    \vspace{0.5cm}
    
    $ \rank \ H_{\Lambda^*}^{B_{1}^{+}}=2, \ and  \ L_{1}=\{0\}$.\\
    
    \vspace{0.5cm}
 
  The algorithm stops with success, the flat extension property is satisfied,\\
     \[ \ker \ H_{\Lambda^*}^{B_{1}^{+}}=\{y-1, x^2-3x+2,xy-1.5y-x+1.5\}\]\\
  and \[ B_{1}=\{1,x-1.5\}. \]

  \end{itemize}
\end{example}

%\subsection{Computing the minimizer points}
\section{Minimizers}\label{sec:5}
In this section we tackle the computation of the minimizer points, once
Algorithm  \ref{algo:flatcrit} stops with \texttt{success} for
$\Lambda^{*} \in \Lc_{E,G}$ optimal  for $f$.
The minimizer points can be computed from the eigenvalues of the
multiplication operators 
$M_{k}: a\in \Ac_{min} \mapsto x_{k} a \in  \Ac_{min}$
for $k=1, \ldots,n$ where $\Ac_{min}= \R[\xx]/I_{min}$ and $I_{min} =
I_{\tilde{\Lambda}}= \Ic(\xi_{1}, \ldots, \xi_{r})$.
At the end of this section we compute the minimizers for our running example.

\begin{prop}\label{multmat}
The matrix of $M_{k}$ in the basis $B$ of $\Ac_{min}$
is $[M_{k}] = (\frac{\Lambda^{*} (x_k \, b_i\, b_j )}{\Lambda^{*}
  (b_i\, b_i)})_{1\le i,j \le r}$.
The operators $M_{k}$, $k=1 \ldots n$ have $r$ common eigenvectors
$\uu_{1}, \ldots, \uu_{r}$ which satisfy $M_{k} \uu_{i} = \xi_{i,k} \uu_{i}$, 
with $\xi_{i,k}$ the $k^{th}$ coordinate of the minimizer point $\xi_{i}= (\xi_{i,1}, \ldots, \xi_{i,n}) \in S$.
\end{prop}
\begin{Proof} 
By Proposition \ref{proposition:algo3.1} and by definition of the
inner-product \eqref{eq:innerprod}, 
$B=\{b_{1}, \ldots, b_{r}\}$ is a basis of $\Ac_{\tilde{\Lambda}}$ and 
$$ 
x_k b_j \equiv \sum_{i = 1}^{r}
        \frac{\Lambda^{*} (x_k \, b_i\, b_j )}{\Lambda^{*} (b_i\,
          b_i)} b_i \mod I_{min},
$$ 
for $j=1\ldots r$, $k=1 \ldots n$. \\
This yields the matrix of the operator $M_{k}$ in the basis $B$: 
$[M_{k}] = (\frac{\Lambda^{*} (x_k \, b_i\, b_j )}{\Lambda^{*} (b_i\,
  b_i)})_{1\le i,j \le r}$.

As the roots of $I_{min}$ are simple, by \citep{em-07-irsea}[Theorem 4.23] the eigenvectors of all
$M_{k}$, $k=1 \ldots n$ are the so-called idempotents 
$\uu_{1}, \ldots, \uu_{r}$ of 
$\Ac_{min}$ and the corresponding eigenvalues are $\xi_{1,k}, \ldots, \xi_{r,k}$.
\end{Proof}

\begin{algorithm2e}[H]\caption{\textsc{Minimizer points}}
\KwIn {$B$ and the output relations from Algorithm \ref{algo:flatcrit}.\\}
\KwOut{ the minimizer points $\xi_{i}= (\xi_{i,1}, \ldots, \xi_{i,n})$, $i=1
  \ldots r$.}
{\tmstrong{Begin}}
\begin{enumerate}
 \item Compute the matrices $[M_{k}] = (\frac{\Lambda^{*} (x_k \, b_i\, b_j )}{\Lambda^{*} (b_i\,
  b_i)})_{1\le i,j \le r}$. 
 \item For a generic choice of $l_{1},\ldots, l_{n}\in \R$, compute the
 eigenvectors $\uu_{1}, \ldots, \uu_{r}$ of $l_{1} [M_{1}] + \cdots +
 l_{n} [M_{n}]$.
 \item Compute $\xi_{i,k}\in \R$ such that $M_{k} \uu_{i} = \xi_{i,k} \uu_{i}$.
\end{enumerate}

{\tmstrong{End}}
% \KwOut{ the minimizer points $\xi_{i}= (\xi_{i,1}, \ldots, \xi_{i,n})$, $i=1
%   \ldots r$.}
\end{algorithm2e}

Now we compute the minimizer points of our running example \ref{runex}.

\begin{example}
 With the basis $B=\{1,x-1.5\}$ and the kernel $\ker \
 H_{\Lambda}^{B^{+}}=\<y-1, x^2-3x+2,xy-1.5y-x+1.5\>$ produced by Algorithm
 \ref{exampledal}, we can compute the multiplication matrices.
  \[M_{x}^{B=\{1,x-1.5\}}=  \left( \begin{array}{cc}
    1.5 & 0.25\\
    1 & 1.5\end{array} \right) \longrightarrow \left\lbrace
  \begin{array}{l}
   x\times 1=1.5 \cdot 1 +1 \cdot (x-1.5)\\
   x\times (x-1.5)=0.25 \cdot 1+1.5 \cdot (x-1.5)
  \end{array}
  \right.  \]
 \[M_{y}^{B=\{1,x-1.5\}}=  \left( \begin{array}{cc}
    1 & 0\\
    0 & 1\end{array} \right) \longrightarrow \left\lbrace
  \begin{array}{l}
   y\times 1=1 \cdot 1 +0 \cdot (x-1.5)\\
   y\times (x-1.5)=0 \cdot 1+1 \cdot (x-1.5)
  \end{array}
  \right.  \]
 
 We take a linear combination of these matrices 
  \[M=M_{x}^{B}+M_{y}^{B}=  \left( \begin{array}{cc}
    2.5 & 0.25\\
    1 & 2.5\end{array} \right) \] 
and compute its eigenvalues
$\lambda_1=2, \lambda_2=3$
and its eigenvectors:
\[M \cdot u_1 =\lambda_1 \cdot u_1 \rightarrow u_1^T=(-0.5,1);\ M \cdot u_2 =\lambda_2 \cdot u_2 \rightarrow u_2^T=(0.5,1)\]
  
 From these eigenvectors, we compute the eigenvalues associated to
 each multiplication matrix $M_{x}^{B}, M_{y}^{B}$. Each computed eigenvalue corresponds 
 to a coordinate of the corresponding minimizer point as we have seen in Proposition \ref{multmat}:
  \[M_{x}^{B} \cdot u_1^T =x_1 \cdot u_1^T \rightarrow x_1=1; \ M_{x}^{B} \cdot u_2^T =x_2 \cdot u_2^T \rightarrow x_2=2\]
  \[M_{y}^{B} \cdot u_1^T =y_1 \cdot u_1^T \rightarrow y_1=1; \ M_{y}^{B} \cdot u_2^T =y_2 \cdot u_2^T \rightarrow y_2=1\]
We recover the minimizer points $(1,1)$ and $(2,1)$.
\end{example}

%\subsection{Minimum with sum of square certification}

%\section{Minimizer border basis algorithm}\label{sec:5}
\section{Main algorithm}\label{sec:7}
In this section we describe the algorithm to compute the 
 infimum of a polynomial on $S$ and the minimizer points when the minimizer ideal is zero-dimensional. %, assuming $f^{*}$ is reached in $\R^{n}$ and
%$I_{min} (f)$ is zero dimensional.
It can be seen as a type of border basis algorithm, in which in the main loop we compute the optimal 
linear form (section \ref{sec:3}), we then check when the minimun is reached (section \ref{sec:4}) and finally
we compute the minimizer points (section \ref{sec:5}). This algorithm is closely connected to the real radical border basis
algorithm presented in \citep{lasserre:hal-00651759}.
%instead of
%  ``minimizing zero'' to generate new elements in the real radical, we
% minimize $f$ to compute generators of the minimizer ideal $I_{min} (f)$.

\begin{algorithm2e}[H]\caption{\textsc{Minimization of $f$ on $S$}\label{algo:bbr}}
\KwIn {A real polynomial function $f$ and a set of constraints
  $\gb\subset \R[\xx]$ with $V_{min}$ non-empty finite.\\}
\KwOut {the minimum $f^{*}= f^{*}_{G_{t},B_{t}}$, the minimizer points $V_{min}=V$, 
  $I_{min}= (K)$ and $B'$ such that $K$ is a border basis for $B'$.}
{\tmstrong{Begin}}
\begin{enumerate}
 \item Take $t=max(\lceil \frac{deg (f)}{2} \rceil,
 d^{0},d^{+})$ where $d^{0}=max_{g^{0} \in \gb^{0}}(\lceil \frac{deg (g^{0})}{2} \rceil), d^{+}=max_{g^{+} \in \gb^{+}}(\lceil \frac{deg (g^{+})}{2} \rceil)$
 \item Compute the graded border basis $F_{2t}$ of $\gb^{0}$ for B in degree $2t$.
 \item Let $B_{t}$ be the set of monomials in $B$ of degree $\le t$.
 \item Let $G_t$ be the set of constraints such that $G^{0}_t=\{m-\pi_{B_t,F_{2t}}(m), \  m \in B_t \cdot B_t \}$ and $G^{+}=\pi_{B_t,F_{2t}}(g^{+})$
 \item $[f^{*}_{G_{t},B_{t}}, \Lambda^{*}]:= \textsc{Optimal Linear Form} (f,B_{t},G_t)$.
 \item %If there is not duality gap 
$[c,B',K]:=\textsc{Decomposition}(\Lambda^{*},B_{t})$ where $c=$\texttt{failed},
$B'=\emptyset, K=\emptyset$  or $c=$\texttt{success}, $B'$ is the basis and $K$ is the set of the relations. 
% \item[] else  go to step 2 with $t:=t+1$
 \item if $c=$\texttt{success} then $V$=\textsc{Minimizer points}$(B',K)$
 \item[] else go to step 2 with $t:=t+1$.
\end{enumerate}
{\tmstrong{End}}

\end{algorithm2e}

\section{Finite convergence}\label{sec:6}

In this section we analyse cases for which an exact relaxation can be constructed.

Our approach to compute the minimizer points relies on the
fact that the border basis relaxation is exact.

By Proposition \ref{lemextend}, the reduced border basis relaxation is exact if and only if the
corresponding full moment matrix relaxation is exact.

Despite the full moment matrix relaxation is not always exact, it is possible to add constraints so that the relaxation becomes exact.

In \citep{Cert}, a general strategy to construct exact SDP relaxation
hierarchies and to compute the minimizer ideal is described. 
It applies to the following problems:

\medskip\noindent{}\textbf{Global optimization.} Consider the case
$n_{1}=n_{2}=0$ with $f^{*}=\inf_{\xx \in \R^{n}} f (\xx)$ reached at a point of $\R^{n}$. 
Taking $G$ such that $G^{0}=\{\frac{\partial f}{\partial_{x_{1}}},
\ldots, \frac{\partial f}{\partial_{x_{n}}} \}$ and $G^{+}=\emptyset$,
the relaxation associated to the sequence
$\Lc_{t,G}$ is exact and
yields $I_{min}$ (see \citep{NDS,Cert}).
If $I_{min}$ is finite then the border basis relaxation yields the minimizer
points and the corresponding border basis.

\medskip\noindent{}\textbf{Regular case.}

We say that $\gb=(g_{1}^{0}, \ldots g_{n_{1}}^{0}$; $g_{1}^{+}, \ldots, g_{n_{2}}^{+})$ is regular if 
for all points $\xx\in \Sc (\gb)$
with $\{j_{1}, \ldots, j_{k}\}=\{ j\in [1, n_{2}] \mid g_{j}^{+} (\xx) =0 \}$, 
the vectors $\nabla g_{1}^{0} (\xx), \ldots, \nabla g_{n_{1}}^{0} (\xx)$,
$\nabla g_{j_{1}}^{+} (\xx), \ldots$, $\nabla g_{j_{k}}^{+} (\xx)$ are linearly independent.

For $\nu =\{j_{1}, \ldots, j_{k}\}\subset [0,n_{2}]$ with $|\nu|\le n-n_{1}$, let 
\begin{eqnarray*}
A_{\nu} &= & [ \nabla f, \nabla g_{1}^{0}, \ldots, \nabla g_{n_{1}}^{0}, \nabla
 g_{j_{1}}^{+}, \ldots, \nabla g_{j_{k}}^{+}]\; \\
\Delta_{\nu} &=& \det (A_{\nu} A_{\nu}^{T}) \\
g_{\nu} &=& \Delta_{\nu} \prod_{j \not \in \nu} g_{j}^{+}.
\end{eqnarray*}
Let $G\subset \R[\xx]$ be the set of constraints such that $G^{0}= \gb^0
\cup \{g_{\nu} \mid \nu \subset [0,n_{2}], |\nu|\le n-n_{1}\}$.
Then the relaxation associated to the preordering sequence $\Lc_{t,G}^{\star}$ is exact and yields
$I_{min}$ (see \citep{Ha-Pham:10,Cert} or \citep{Nie11} for
$\C$-regularity and constraints $G^{0}$ that involve minors of $A_{\nu}$).

If $I_{min}$ is non-empty and finite then the border basis relaxation
\eqref{relaxation:bb} yields the points $V_{min}$ and the border basis of $I_{min}$.

\medskip\noindent{}\textbf{Boundary Hessian Conditions.} If $f$ and $\gb$ satisfies the so-called Boundary Hessian Conditions %at each zero of $f-f^*$ on $S$
then $f-f^* \in \Qc_{t,g}$ and the relaxation associated to $\Lc_{t,\gb}$ is exact and yields $I_{min}$ (see \citep{Marshall}).
If moreover $I_{min}$ is finite then the border basis relaxation yields the
points $V_{min}$ and the corresponding border basis of $I_{min}$.

\medskip\noindent{}\textbf{$\gb^{+}$-radical computation.}
If we optimize $f=0$ on the set $S = \Sc (\gb)$, then all the points of $S$
are minimizer points, $V_{min} = S$ and by the Positivstellensatz, $I_{min}$
is equal to 
$$
\sqrt[\gb^{+}]{\gb^{0}}=\{p \in \R[\xx] \mid \exists m \in \N
\ s.t. \ p^{2m} + q=0, q \in \Pc_{\R[\xx],\gb} \}.
$$
Here again, the preordering sequence $\Lc_{{t},\gb}^{\star}$ is exact. If we assume that $S=\Sc (\gb)$ is finite, then the corresponding border basis
relaxation yields the points of $S$ and the generators of 
$\sqrt[\gb^{+}]{\gb^{0}}$.
See also \citep{LLR08b, lasserre:hal-00651759} for zero dimensional
real radical computation and \citep{Zhi}.

%\subsection{Experimentation}\label{sec:experiment}
\section{Performance}\label{sec:experiment}
In this section, we analyse the practical behavior of Algorithm \ref{algo:bbr}. In all the examples the minimizer
ideal is zero-dimensional hence our algorithm stops in a finite
number of steps and yields the minimizer points and generators of the
minimizer ideal.

The implementation of the previous algorithm has been performed using
the {\sc borderbasix}\footnote{http://www-sop.inria.fr/teams/galaad/software/bbx/} package of the {\sc Mathemagix}\footnote{www.mathemagix.org} software,
which provides a {\tt C++} implementation of the border basis algorithm of
\citep{BMPhT12}.

For the computation of border basis, we use  a choice function that is tolerant
to numerical inestability i.e. a choice function that chooses as
leading monomial a monomial whose coefficient is maximal among the
choosable monomials as described in \citep{BMPhT08}.
% makes the border
% basis computation stable with respect to numerical perturbations.

The Semi-Definite Programming problems are solved using
\textsc{sdpa},
\textsc{sdpa-gmp}\footnote{http://sdpa.sourceforge.net}, \textsc{csdp}
and \textsc{mosek}\footnote{http://www.mosek.com} software. For the link with \textsc{sdpa},\textsc{csdp} and \textsc{sdpa-gmp} we use a file interface. In the case of \textsc{mosek}, we use the distributed library. 
% The Semi-Definite Programming problems are solved using
%  \textsc{mosek}\footnote{http://www.mosek.com}software. For the
% link with \textsc{mosek}, we use the distributed library.

Once we have computed the moment matrix, we call the Decomposition Algorithm which is available 
in the {\sc borderbasix} package.

The minimizer points are computed from the eigenvalues of the multiplication
matrices. This is performed using Lapack routines.

% After computing the basis and the border basis relations of the
% minimizer ideal, the roots are
% obtained using the numerical routines described in \citep{SGPhT09}.

Experiments are made on an Intel Core i5 2.40GHz. %with 8Gb of RAM.

In Table \ref{tabla}, we compare our algorithm \ref{algo:bbr} (\textit{bbr})
with the full moment matrix relaxation algorithm (\textit{fmr}) inside the same environment.
This latter (implemented by ourselves in \texttt{C++} inside the \textsc{borderbasix} package) reproduces the algorithm described in \citep{Lasserre:book}, which
is also implemented in the package \textsc{gloptipoly} of \textsc{matlab} developed by
D. Henrion and J.B. Lasserre.
In this table, we record the problem name or the source of the problem,
the number of decision variables (v), the number of inequality and
equality constraints, we mark in parenthesis the number of equality constraints (c), the maximum degree of 
the constraints and of the polynomial to minimize (d), the number of minimizer points (sol). For the
two algorithms \textit{bbr} and \textit{fmr} we report the total CPU
time in seconds using \textsc{mosek} (t),
the order of the relaxation (o), the number of parameters of the SDP problem
(p) and the size of the moment matrices (s). 
The first part of the table contains examples of positive polynomials,
which are not sums of squares. New equality constraints are added following \citep{Cert} to
compute the minimizer points in the examples marked with $\diamond$.
% The fourth part of the table contains examples where the real radical
% $\sqrt[\gb^{+}]{\gb^{0}}$ is computed.
When there are equality constraints, the border basis computation
reduces the size of the moment matrices, as well as the 
localization matrices associated to the inequalities. 
This speeds up the SDP computation as shown the examples Ex 1.4.8, Ex 2.1.8, Ex 2.1.9 and simplex.
In the case where there are only inequalities, the size of the moment
matrices and number of parameters do not change but once the optimal linear form is computed
using the SDP solver \textsc{mosek},
the \textsc{Decomposition} algorithm which computes the minimizers is
more efficient and quicker than the reconstruction algorithm used 
in the full moment matrix relaxation approach.
The performance is not the only issue: numerical problems can also occur
due to the bigger size of the moment matrices in the flat extension
test and the reconstruction of minimizer points. Such examples where the
\textit{fmr} algorithm fails are marked with \textsc{*}. In these three
problems, there is not a big enough gap between the singular values to
determine correctly the numerical rank and the flat extension property cannot
be verified. 
The examples that \textsc{gloptipoly} cannot treat due to the high number of
variables \citep{Lasserre:book} are marked with \textsc{**}. We can treat three of this examples (with \textit{fmr}) 
because as we said \textit{fmr} is implemented in \texttt{C++} so it is more
efficient than \textsc{glotipoly}, which is implemented inside \textsc{matlab}.
We cannot treat the example 2.1.8 with the \textit{fmr} algorithm due
to the large number of parameters.
%We remark that these examples can be solved with \textit{fmr} showing
% a better performance in C++ in comparison with matlab.

These experiments show that when the size of the SDP problems becomes significant,
most of the time is spent during the SDP computation and the border
 basis time and reconstruction time are negligible. The use of
 \textsc{mosek} software provides a speed-up factor of 1.5 to 5
 compared to the \textsc{sdpa} software for small examples (such as
 Robinson, Moztkin, Ex 3, Ex 5, Ex 2.1.1, Ex 2.1.2, Ex 2.1.4 and Ex
 2.1.6). For large examples (such as Ex 2.1.3, Ex 2.1.7, Ex 2.1.8 and
 simplex) the improvement factors are between 10-30 times. 
 These improvements are due to the new fast Cholesky decomposition inside of \textsc{mosek} software \footnote{http://www.mosek.com}.
 %We also show that the use of mosek software
% reduces the time between 50 \% and 80 \%.
In all the examples, the new border basis relaxation algorithm
outperforms the full moment matrix relaxation method.

In Table \ref{tabla2}, we apply our algorithm \textit{bbr} to find the best rank-1 and rank-2 tensor approximation for
symmetric and non symmetric tensors on examples from 
\citep{Nietensors} and \citep{Structmatrix}. For best rank-1
approximation problems with several minimizers (which is the case when
there are symmetries), the method proposed in \citep{Nietensors}
cannot certify the result and uses a local method to converge to a
local extrema. We apply the global border basis relaxation algorithm to find all the 
minimizers for the best rank 1 approximation problem.

The last example in Table \ref{tabla2} is
a best rank-2 tensor approximation example from the paper
\citep{Structmatrix}. The eight solutions come from the symmetries due
to the invariance of the solution set by permutation and negation of the factors.

%\begin{changemargin}{cm}{2cm} 
\begin{table*}
\ \ \ \ \ \ \ \ 
\begin{changemargin}{-1.5cm}{2cm} 
\begin{tabular}{|c||ccc||c||c|c|c|c||c|c|c|c|}\hline
 problem & v & c & d & sol & $o_{bbr}$& $p_{bbr}$&$s_{bbr}$&$t_{bbr}$&$o_{fmr}$& $p_{fmr}$&$s_{fmr}$&$t_{fmr}$ \\ \hline
 $\diamond$ Robinson              & 2  & 0  & 6 & 8 &4 &20 &15 &0.07&7&119&36&*\\
 $\diamond$ Motzkin               & 2  & 0  & 6 & 4 &4 &25 &15 &0.060&9&189&55&*\\
 $\diamond$  Motzkin perturbed    & 3  & 1  & 6 & 1 &5 &127&35 &0.18&5&286&56&8.01\\\hline
 $\diamond$ L'01, Ex. 1           & 2  & 0  & 4 & 1 &2 &8  &6  &0.020&2&14&6&0.035\\
 $\diamond$ L'01, Ex. 2           & 2  & 0  & 4 & 1 &2 &8  &6  &0.020&2&14&6&0.026\\
 $\diamond$ L'01, Ex. 3           & 2  & 0  & 6 & 4 &4 &25 &15 &0.057&8&152&45&*\\
  L'01, Ex. 5                     & 2  & 3  & 2 & 3 &2 &14 &6  &0.032&2&14&6&0.045\\ \hline
 % \citep{Handbooktest}, Ex. 4.1.1 & 1  & 2  & 6 & 1  &0.027&3&6&4&0.048&3&6&4 \\
 %\citep{Handbooktest}, Ex. 4.1.3 & 1  & 2  & 5 & 1  &0.024&3&6&4&0.038&3&6&4 \\
 F, Ex. 4.1.4                     & 1  & 2  & 4 & 2 &2 &4  &3  &0.016&2&4&3&0.023\\
 %\citep{Handbooktest}, Ex. 4.1.5 & 2  & 2  & 6 & 1  &0.060&3&27&10&*&3&27&10\\
 F, Ex. 4.1.6                     & 1  & 2  & 6 & 2 &3 &6&4   &0.018&3&6&4&0.020\\
 F, Ex. 4.1.7                     & 1  & 2  & 4 & 1 &2 &4&3 &0.017&2&4&3&0.020\\
 F, Ex. 4.1.8                     & 2  & 5(1)  & 4 & 1 &2 &13&6 &0.021&2&14&6&0.12\\
%  F, Ex. 4.1.9                     & 2  & 6  & 4 & 1 &4 &44&15  &0.11&4&44&15&0.29\\ 
 F, Ex. 2.1.1                     & 5  & 11 & 2 & 1 &3 &461&56 &3.10&3&461&56&3.12\\
 F, Ex. 2.1.2                     & 6  & 13 & 2 & 1 &2 &209&26 &0.32&2&209&26&0.36\\
 F, Ex. 2.1.3                     & 13 & 35 & 2 & 1 &2 &2379&78 &19.68&2&2379&78&25.60\\
 F, Ex. 2.1.4                     & 6  & 15 & 2 & 1 &2 &209&26  &0.30&2&209&26&0.33\\
 F, Ex. 2.1.5                     & 10 & 31 & 2 & 1 &2 &1000&66 &9.15&2&1000&66&9.7\\
 F, Ex. 2.1.6                     & 10 & 25 & 2 & 1 &2 &1000&66 &3.6&2&1000&66&4.17\\
 **F, Ex. 2.1.7(1)                & 20 & 30 & 2 & 1 &2 &10625&231 &730.24&2&10625&231&1089.31 \\ 
 ** F, Ex. 2.1.7(5)               & 20 & 30 & 2 & 1 &2&10625&231 &747.94&2&10625&231&1125.27 \\ 
  ** F, Ex. 2.1.8                 & 24 & 58(10) & 2 & 1 &2&3875&136  &311.54&2&20474&325&>14h \\ 
 F, Ex. 2.1.9                     & 10 & 11(1) & 2 & 1 &2&714&44  &0.62&2&1000&55&1.67\\ \hline
 % ** \citep{Handbooktest}, Ex. 2.1.10 & 20 & 30 & 2 & 1  &-&2&10625&231&-&2&10625&231 \\ 
%  F, Ex. 3.1.3                     & 6  & 16 & 2 & 1 & 2&209&26  &0.61&2&209&26&1.42\\ \hline
 %\citep{Handbooktest}, Ex. 4.1.2 & 1  & 2  &50 & 1  &0.52&25&50&25&2.38&26&52&26\\
% %  L'09 cbms1		          & 3  & 3  & 3 & 5 &3&26&17 &0.14&3&83&20&0.20\\
% %  L'09 rediff3		          & 3  & 3  & 2 & 2 &2&7&7 &0.06&2&35&10&0.09\\
% %  L'09 quadfor2                    & 4  & 12 & 4 & 2 &3&48&19 &0.45&3&210&35&0.75\\ \hline
  
 ** simplex                       & 15 & 16(1) & 2 & 1 &2&3059&120 &15.30&2&3875&136&47.50\\ 
\hline 
\end{tabular}

\center{\caption{Examples from F-\citep{Handbooktest}),  L'01-\citep{Las01}.\label{tabla}}}%L'09-\citep{Lasserre:book},
\end{changemargin}
\end{table*}

\begin{table*}
%\ \ \ \ \ \ \ \ 
\begin{tabular}{|c||ccc||c||c|c|c|c|}\hline
 problem & v & c & d & sol &$o_{bbr}$& $p_{bbr}$&$s_{bbr}$&$t_{bbr+msk}$\\ \hline
 \citep{Nietensors} Ex. 3.1 & 2& 1& 3& 1&2&8&5&0.028\\
 \citep{Nietensors} Ex. 3.2 & 3& 1& 3& 1&2&24&9&0.025\\
 \citep{Nietensors} Ex. 3.3 & 3& 1& 3& 1&2&24&9&0.035\\
 \citep{Nietensors} Ex. 3.4 & 4& 1& 4& 2&2&24&9&0.097\\
 \citep{Nietensors} Ex. 3.5 & 5& 1& 3& 1&2&104&20&0.078\\
 \citep{Nietensors} Ex. 3.6 & 5& 1& 4& 2&4&824&105&15.39\\
 \citep{Nietensors} Ex. 3.8 & 3& 1& 6& 4&3&48&16&1.14\\
 \citep{Nietensors} Ex. 3.11 & 8& 4& 4& 8&3&84&25&0.17\\
 \citep{Nietensors} Ex. 3.12 & 9& 3& 3& 4&2&552&52&1.55\\
 \citep{Nietensors} Ex. 3.13 & 9& 3& 3& 12&3&3023&190&223.27\\ \hline 
 \citep{Structmatrix} Ex. 4.2 & 6&0&8&4&8&2340&210&59.38\\ \hline
\end{tabular}

\center{\caption{Best rank-1 and rank-2 approximation tensors \label{tabla2}}}
\end{table*}

%\medskip\noindent{}\textbf{Acknowledgements.} We would like to thank Philippe Trebuchet and Matthieu Dien for their development in the {\sc borderbasix} package.
\section*{Acknowledgements}
 We would like to thank Philippe Trebuchet and Matthieu Dien for their development in the {\sc borderbasix} package.

%\section*{References}

%\bibliographystyle{plain}

%\bibliography{paper}

\appendix
\section{Results of Best rank-1 approximation tensors}
\noindent{}\textit{Example 3.1}: Consider the tensor $\Fc \in S^3(\R^2)$ with entries\\ 
 $\Fc_{111}=1.5578, \Fc_{222}=1.1226, \Fc_{112}=-2.443, \Fc_{221}=-1.0982$\\
 We get the rank-1 tensor $\lambda \cdot u^{\otimes 3}$ with:\\
 $\lambda =3.11551, \ u=(0.926433,-0.376457)$  and  $\mid \mid \Fc -\lambda \cdot u^{\otimes 3}\mid \mid=3.9333.$\\

\noindent{}\textit{Example 3.2}: Consider the tensor  $\Fc \in S^3(\R^3)$ with entries \\
 $\Fc_{111}=-0.1281, \Fc_{112}=0.0516, \Fc_{113}=-0.0954, \Fc_{122}=-0.1958, \Fc_{123}=-0.1790,$\\
 $\Fc_{133}=-0.2676, \Fc_{222}=0.3251, \Fc_{223}=0.2513, \Fc_{233}=0.1773, \Fc_{333}=0.0338$\\
 We get the rank-1 tensor $\lambda \cdot u^{\otimes 3}$ with:\\
 $\lambda =0.87298, \ u=(-0.392192,0.7248,0.5664)$ and  
 $\mid \mid \Fc -\lambda \cdot u^{\otimes 3}\mid \mid=0.4498$.\\
 
\noindent{}\textit{Example 3.3}: Consider the tensor  $\Fc \in S^3(\R^3)$ with entries \\
 $\Fc_{111}=0.0517, \Fc_{112}=0.3579, \Fc_{113}=0.5298, \Fc_{122}=0.7544, \Fc_{123}=0.2156,$\\
 $\Fc_{133}=0.3612, \Fc_{222}=0.3943, \Fc_{223}=0.0146, \Fc_{233}=0.6718, \Fc_{333}=0.9723$\\
  We get the rank-1 tensor $\lambda \cdot u^{\otimes 3}$ with:\\
 $\lambda =2.11102, \ u=(0.52048,0.511264,0.683891)$ and  
 $\mid \mid \Fc -\lambda \cdot u^{\otimes 3}\mid \mid=1.2672$.\\

\noindent{}\textit{Example 3.4}: Consider the tensor  $\Fc \in S^4(\R^3)$ with entries \\
 $\Fc_{1111}=0.2883, \Fc_{1112}=-0.0031, \Fc_{1113}=0.1973, \Fc_{1122}=-0.2458, \Fc_{1123}=-0.2939,$\\
 $\Fc_{1133}=0.3847, \Fc_{1222}=0.2972, \Fc_{1223}=0.1862, \Fc_{1233}=0.0919, \Fc_{1333}=-0.3619$\\
 $\Fc_{2222}=0.1241, \Fc_{2223}=-0.3420, \Fc_{2233}=0.2127, \Fc_{2333}=0.2727, \Fc_{3333}=-0.3054$\\
  We get the rank-1 tensor $\lambda \cdot u_i^{\otimes 3}$ with:\\
 $\lambda =-1.0960, \ u_1=(-0.59148,0.7467,0.3042); u_2=(0.59148,-0.7467,-0.3042)$ and  
 $\mid \mid \Fc -\lambda \cdot u_i^{\otimes 4}\mid \mid=1.9683$.\\
 
\noindent{}\textit{Example 3.5}: Consider the tensor  $\Fc \in S^3(\R^5)$ with entries \\
 $\Fc_{i_1,i_2,i_3}=\frac{(-1)^{i_1}}{i_1}+\frac{(-1)^{i_2}}{i_2}+\frac{(-1)^{i_3}}{i_3}$\\
  We get the rank-1 tensor $\lambda \cdot u^{\otimes 3}$ with:\\
 $\lambda =9.9776, \ u=(-0.7313,-0.1375,-0.46737,-0.23649,-0.4146)$ and \\
 $\mid \mid \Fc -\lambda \cdot u^{\otimes 3}\mid \mid=5.3498$.\\
 
\noindent{}\textit{Example 3.6}: Consider the tensor  $\Fc \in S^4(\R^5)$ with entries \\
 $\Fc_{i_1,i_2,i_3,i_4}=arctan((-1)^{i_1}\frac{i_1}{5})+arctan((-1)^{i_2}\frac{i_2}{5})+arctan((-1)^{i_3}\frac{i_3}{5})+arctan((-1)^{i_4}\frac{i_4}{5})$\\
  We get the rank-1 tensor $\lambda \cdot u^{\otimes 4}$ with:\\
 $\lambda =-23.56525, \ u_1=(0.4398,0.2383,0.5604,0.1354,0.6459);\\ \ u_2=(-0.4398,-0.2383,-0.5604,-0.1354,-0.6459)$ and \\
 $\mid \mid \Fc -\lambda \cdot u_i^{\otimes 4}\mid \mid=16.8501$.\\

\noindent{}\textit{Example 3.8}: Consider the tensor  $\Fc \in S^6(\R^3)$ with entries \\
 $\Fc_{111111}=2, \Fc_{111122}=1/3, \Fc_{111133}=2/5, \Fc_{112222}=1/3, \Fc_{112233}=1/6,$\\
 $\Fc_{113333}=2/5, \Fc_{222222}=2, \Fc_{222233}=2/5, \Fc_{223333}=2/5, \Fc_{333333}=1$\\
  We get the rank-1 tensor $\lambda \cdot u_i^{\otimes 6}$ with:\\
 $\lambda =2, \ u_1=(1,0,0); \ u_2=(-1,0,0); \ u_3=(0,1,0); \ u_4=(0,-1,0)$ and \\
 $\mid \mid \Fc -\lambda \cdot u_i^{\otimes 6}\mid \mid=20.59$.\\
 
\noindent{}\textit{Example 3.11}: Consider the tensor  $\Fc \in \R^{2 \times 2 \times 2 \times 2}$ with entries\\
  $\Fc_{1111}=25.1, \Fc_{1212}=25.6, \Fc_{2121}=24.8, \Fc_{2222}=23$\\
 We get the rank-1 tensor $\lambda \cdot u_i^1\otimes u_i^2 \otimes u_i^3 \otimes u_i^4$ with:\\
 $\lambda =25.6, \ u_1^1=(1,0), u_1^2=(0,1),  u_1^3=(1,0), u_1^4=(0,1);\\
 \ u_2^1=(-1,0), u_2^2=(0,-1),  u_2^3=(-1,0), u_2^4=(0,-1);\\
 \ u_3^1=(-1,0), u_3^2=(0,-1),  u_3^3=(1,0), u_3^4=(0,1);\\
 \ u_4^1=(1,0), u_4^2=(0,1),  u_4^3=(-1,0), u_4^4=(0,-1);\\
 \ u_5^1=(-1,0), u_5^2=(0,1),  u_5^3=(-1,0), u_5^4=(0,1);\\
 \ u_6^1=(1,0), u_6^2=(0,-1),  u_6^3=(1,0), u_6^4=(0,-1);\\
 \ u_7^1=(1,0), u_7^2=(0,-1),  u_7^3=(-1,0), u_7^4=(0,1);\\
 \ u_8^1=(-1,0), u_8^2=(0,1),  u_8^3=(1,0), u_8^4=(0,-1)$.\\
 The distance between $\Fc$ and one of these solutions is
 $\mid \mid \Fc -\lambda \cdot u_i^1\otimes u_i^2 \otimes u_i^3 \otimes u_i^4 \mid \mid=42.1195$.\\
 
\noindent{}\textit{Example 3.12}: Consider the tensor  $\Fc \in \R^{3 \times 3 \times 3}$ with entries\\
 $\Fc_{111}=0.4333, \Fc_{121}=0.4278, \Fc_{131}=0.4140, \Fc_{211}=0.8154, \Fc_{221}=0.0199,$\\
 $\Fc_{231}=0.5598, \Fc_{311}=0.0643, \Fc_{321}=0.3815, \Fc_{331}=0.8834, \Fc_{112}=0.4866,$\\
 $\Fc_{122}=0.8087, \Fc_{132}=0.2073, \Fc_{212}=0.7641, \Fc_{222}=0.9924, \Fc_{232}=0.8752,$\\
 $\Fc_{312}=0.6708, \Fc_{322}=0.8296, \Fc_{332}=0.1325, \Fc_{113}=0.3871, \Fc_{123}=0.0769,$\\
 $\Fc_{133}=0.3151, \Fc_{213}=0.1355, \Fc_{223}=0.7727, \Fc_{233}=0.4089, \Fc_{313}=0.9715,$\\
 $\Fc_{323}=0.7726, \Fc_{333}=0.5526$\\
  We get the rank-1 tensor $\lambda \cdot u_i^1\otimes u_i^2 \otimes u_i^3$ with:\\
 $\lambda =2.8166, \ u_1^1=(0.4279,0.6556,0.62209), u_1^2=(0.5705,0.6466,0.5063), u_1^3=(0.4500,0.7093,0.5425);\\
 \ u_2^1=(0.4279,0.6556,0.62209), u_2^2=(-0.5705,-0.6466,-0.5063), u_2^3=(-0.4500,-0.7093,-0.5425);\\
 \ u_3^1=(-0.4279,-0.6556,-0.62209), u_3^2=(0.5705,0.6466,0.5063), u_3^3=(-0.4500,-0.7093,-0.5425);\\
 \ u_4^1=(-0.4279,-0.6556,-0.62209), u_4^2=(-0.5705,-0.6466,-0.5063), u_4^3=(0.4500,0.7093,0.5425),$ \\ 
 The distance between $\Fc$ and one of these solutions is $\mid \mid \Fc -\lambda \cdot u_i^1\otimes u_i^2 \otimes u_i^3  \mid \mid=1.3510$.\\
 
\noindent{}\textit{Example 3.13}: Consider the tensor  $\Fc \in \R^{3 \times 3 \times 3}$ with entries\\
 $\Fc_{111}=0.0072, \Fc_{121}=-0.4413, \Fc_{131}=0.1941, \Fc_{211}=-04413, \Fc_{221}=0.0940,$\\
 $\Fc_{231}=0.5901, \Fc_{311}=0.1941, \Fc_{321}=-0.4099, \Fc_{331}=-0.1012, \Fc_{112}=-0.4413,$\\
 $\Fc_{122}=0.0940, \Fc_{132}=-0.4099, \Fc_{212}=0.0940, \Fc_{222}=0.2183, \Fc_{232}=0.2950,$\\
 $\Fc_{312}=0.5901, \Fc_{322}=0.2950, \Fc_{332}=0.2229, \Fc_{113}=0.1941, \Fc_{123}=0.5901,$\\
 $\Fc_{133}=-01012, \Fc_{213}=-0.4099, \Fc_{223}=0.2950, \Fc_{233}=0.2229, \Fc_{313}=-0.1012,$\\
 $\Fc_{323}=0.2229, \Fc_{333}=-0.4891$\\
  We get the rank-1 tensor $\lambda \cdot u_i^1\otimes u_i^2 \otimes
  u_i^3$ with 
 $\lambda =1.000$ and the 12 solutions\\
 $\ u_1^1=(0.7955,0.2491,0.5524), u_1^2=(-0.0050,0.9142,-0.4051), u_1^3=(-0.6060,0.3195,0.7285);\\
 \ u_2^1=(-0.0050,0.9142,-0.4051), u_2^2=(-0.6060,0.3195,0.7285),u_2^3=(0.7955,0.2491,0.5524);\\
 \ u_3^1=(-0.6060,0.3195,0.7285),u_3^2=(0.7955,0.2491,0.5524), u_3^2=(-0.0050,0.9142,-0.4051);\\
  \ u_4^1=(0.7955,0.2491,0.5524), u_4^2=(0.0050,-0.9142,0.4051), u_4^3=(0.6060,-0.3195,-0.7285);\\
  \ u_5^1=(0.6060,-0.3195,-0.7285),u_5^2=(0.7955,0.2491,0.5524), u_5^2=(0.0050,-0.9142,0.4051);\\
  \ u_6^1=(-0.6060,0.3195,0.7285),u_6^2=(-0.7955,-0.2491,-0.5524), u_6^2=(0.0050,-0.9142,0.4051);\\
  \ u_7^1=(0.6060,-0.3195,-0.7285),u_7^2=(-0.7955,-0.2491,-0.5524), u_7^2=(-0.0050,0.9142,-0.4051);\\
  \ u_8^1=(-0.7955,-0.2491,-0.5524), u_8^2=(-0.0050,0.9142,-0.4051), u_8^3=(0.6060,-0.3195,-0.7285);\\
  \ u_9^1=(-0.7955,-0.2491,-0.5524), u_9^2=(0.0050,-0.9142,0.4051), u_9^3=(-0.6060,0.3195,0.7285);\\
  \ u_{10}^1=(-0.0050,0.9142,-0.4051), u_{10}^2=(0.6060,-0.3195,-0.7285),u_{10}^3=(-0.7955,-0.2491,-0.5524);\\
  \ u_{11}^1=(0.0050,-0.9142,0.4051), u_{11}^2=(0.6060,-0.3195,-0.7285),u_{11}^3=(0.7955,0.2491,0.5524);\\
  \ u_{12}^1=(0.0050,-0.9142,0.4051), u_{12}^2=(-0.6060,0.3195,0.7285),u_{12}^3=(-0.7955,-0.2491,-0.5524).$\\
 The distance between $\Fc$ and one of these solutions is
 $\mid \mid \Fc -\lambda \cdot u_i^1\otimes u_i^2 \otimes u_i^3  \mid \mid=1.4143$.\\
 
\noindent{}\textit{Example 4.2}: Consider the tensor  $\Fc \in S^4(\R^3)$ with entries \\
 $\Fc_{1111}=0.1023, \Fc_{1112}=-0.002, \Fc_{1113}=0.0581, \Fc_{1122}=0.0039, \Fc_{1123}=-0.00032569,$\\
 $\Fc_{1133}=0.0407, \Fc_{1222}=0.0107, \Fc_{1223}=-0.0012, \Fc_{1233}=-0.0011, \Fc_{1333}=0.0196,$\\
 $\Fc_{2222}=0.0197, \Fc_{2223}=-0.0029, \Fc_{2233}=-0.00017418, \Fc_{2333}=-0.0021,$\\
 $\Fc_{3333}=0.1869$\\
  We get the rank-2 tensor
  $\tilde{\Fc}(s,t,u)=(as+bt+cu)^4+(ds+et+fu)^4$ with the 8 solutions:\\
 $s_1=(a,b,c,d,e,f)=(0.01877,0.006239,-0.6434,-0.5592,0.008797,-0.3522);$\\
 $s_2=(-0.01877,-0.006239,0.6434,0.5592,-0.008797,0.3522);$\\
 $s_3=(0.01877,0.006239,-0.6434,0.5592,-0.008797,0.3522);$\\
 $s_4=(-0.01877,-0.006239,0.6434,-0.5592,0.008797,-0.3522);$ \\ 
 $s_5=(-0.5592,0.008797,-0.3522,0.01877,0.006239,-0.6434);$\\
 $s_6=(0.5592,-0.008797,0.3522,-0.01877,-0.006239,0.6434);$\\
 $s_7=(-0.5592,0.008797,-0.3522,-0.01877,-0.006239,0.6434);$\\
 $s_8=(0.5592,-0.008797,0.3522,0.01877,0.006239,-0.6434);$. \\
 The distance between $\Fc$ and one of these solutions is
 $\mid \mid \Fc -\tilde{\Fc} \mid \mid=0.00108483$.\\
The other possible real rank-2 approximations $\tilde{\Fc}(s,t,u)=\pm(as+bt+cu)^4\pm(ds+et+fu)^4$
yield solutions which are not as close to $\Fc$ as these solutions.
\end{document}